\documentclass{article}
\usepackage{latexsym,amsmath,amssymb,amscd,bm,fancyheadings, pifont }
\usepackage{mathrsfs}
\usepackage[dvips]{color}
\usepackage[dvips]{graphicx}
\begin{document}

\title{On the Essential Spectrum of Adiabatic Stellar Oscillations \\
Dedicated to Professor Shih-Hsien Yu to celebrate his sixtieth birthday}
\author{Tetu Makino \footnote{Professor Emeritus at Yamaguchi University, Japan; E-mail: makino@yamaguchi-u.ac.jp}}
\date{\today}
\maketitle

\newtheorem{Lemma}{Lemma}
\newtheorem{Proposition}{Proposition}
\newtheorem{Theorem}{Theorem}
\newtheorem{Definition}{Definition}
\newtheorem{Remark}{Remark}
\newtheorem{Corollary}{Corollary}
\newtheorem{Notation}{Notation}
\newtheorem{Assumption}{Assumption}
\newtheorem{Question}{Question}
\newtheorem{Supposition}{Supposition}

\numberwithin{equation}{section}

\begin{abstract}
The generator $\bm{L}$ of the linearized evolution equation of adiabatic oscillations of a gaseous star, ELASO, is a second order integro-differential operator and is realized as a self-adjoint operator in the Hilbert space of square integrable unknown functions with weight, which is the density distribution of the compactly supported background. Eigenvalues and eigenfunctions of the operator $\bm{L}$ have been investigated in practical point of view of eigenmode expansion of oscillations. But it should be examined whether continuous spectra are absent in the spectrum of $\bm{L}$ or not. In order to discuss this question, the existence of essential spectra in a closely related evolution problem is established.



\end{abstract}

\section{Introduction}
We consider the equation of linearized adiabatic stellar oscillations, {\bf ELASO }:

\begin{equation}
\frac{\partial^2\bm{u}}{\partial t^2}+\mathcal{B} \frac{\partial \bm{u}}{\partial t}+\mathcal{L}\bm{u}=0, \quad t \geq 0, \bm{x} \in \mathfrak{R}_b, \label{ELASO}
\end{equation}
where the unknown is $\bm{u}=\bm{u}(t,\bm{x}) \in \mathbb{R}^3$,
$\mathfrak{R}_b=\{ \bm{x} | \rho_b(\bm{x}) >0\}$ is a bounded domain in $\mathbb{R}^3$.
 and
\begin{subequations}
\begin{align}
\mathcal{B}\bm{v}&=2\Omega
\begin{bmatrix}
-v^2 \\
\\
v^1 \\
\\
0
\end{bmatrix}
, \\
\mathcal{L}\bm{u}&=\mathcal{L}_0\bm{u}+4\pi\mathsf{G}\mathcal{L}_1\bm{u}, \\
\mathcal{L}_0\bm{u}&=\frac{1}{\rho_b}\nabla\delta P -\frac{\nabla P_b}{\rho_b^2}\delta \rho, \\
&\delta\rho=-\mathrm{div}(\rho_b \bm{u}), \quad
\delta P=\frac{\gamma P_b}{\rho_b}\delta\rho +\gamma P_b(\bm{u}|\mathfrak{a}_b), \\
&\mathfrak{a}_b=-\frac{1}{\gamma \mathsf{C}_V}\nabla S_b = \frac{\nabla \rho_b}{\rho_b}-
\frac{\nabla P_b}{\gamma P_b}, \\
\mathcal{L}_1\bm{u}&=\nabla \mathcal{K}[\delta \rho], \\
&\mathcal{K}[g](\bm{x})=\frac{1}{4\pi}\int_{\mathfrak{R}_b}\frac{g(\bm{x}')}{\|\bm{x}-\bm{x}'\|}d\bm{x}'.
\end{align}
\end{subequations}

Here $\mathsf{G}, \gamma, \mathsf{C}_V$ are positive constants, $1<\gamma <2$, and
$(\rho, S, \bm{v})=(\rho_b(\bm{x}), S_b(\bm{x}), \bm{0})$ is a stationary solution of the Euler-Poisson equations in the rotating co-ordinates with a constant angular velocity
$\Omega$:
\begin{subequations}
\begin{align}
&\frac{D\rho}{D t}+\rho\mathrm{div}\bm{v}=0, \label{EPa}\\
&\rho\Big[\frac{D \bm{v}}{Dt} +\bm{\Omega}\times\bm{v}+
\bm{\Omega}\times(\bm{\Omega}\times \bm{x})\Big]+
\mathrm{grad} P+\rho\mathrm{grad}\Phi=0, \\
&\rho\frac{DS}{Dt}=0, \\
&\Phi(t,\dot)=4\pi\mathsf{G}\mathcal{K}[\rho(t,\cdot)], \label{EPd}
\end{align}
\end{subequations}
where $\displaystyle \bm{\Omega}=\Omega\frac{\partial}{\partial x^3},$
$\displaystyle \frac{D}{Dt}=\frac{\partial}{\partial t}+\sum v^k\frac{\partial}{\partial x^k}$, under the equation of state
\begin{equation}
P=\rho^{\gamma}e^{\frac{S}{\mathsf{C}_V}}. \label{EOS}
\end{equation}
The equations \eqref{EPa} - \eqref{EPd}, \eqref{EOS} govern the adiabatic inviscid interior motion of a gaseous star, where $\rho \geq 0$ is the density, $P$ the pressure, $S$ the specific entropy, $\bm{v}$ the velocity field, and $\Phi$ is the gravitational potential.

We assume that $\mathfrak{R}_b=\{ \bm{x} | \rho_b(\bm{x})>0 \}$ is a bounded domain of class $C^{3,\alpha}$, and 
$\rho_b^{\gamma-1}, S_b \in C^{\infty}(\mathfrak{R}_b) \cap C^{3,\alpha}(\mathfrak{R}_b \cup \partial\mathfrak{R}_b)$, $\alpha$ being a positive number such that $0<\alpha< \Big(\frac{1}{\gamma-1}-1\Big)\wedge 1$, 
\begin{equation}
\inf_{0<r<r_0}\Big(-\frac{1}{r}\frac{\partial \rho_b}{\partial r}\Big) >0\quad\mbox{for}\quad
0<r_0 \ll 1,
\end{equation}
where $r=|\bm{x}|$ while we are looking
$$\rho_b=\rho_b(r\sin\vartheta\cos\phi, r\sin\vartheta\sin\phi, r\cos\vartheta),$$
and
\begin{equation}
-\infty< \frac{\partial {c_b^2}}{\partial \bm{n}} < 0 \quad\mbox{on}\quad \partial\mathfrak{R}_b,
\end{equation}
where $\bm{n}$ is the outer normal vector at the boundary point, and
\begin{equation}
{c_b}=\sqrt{ \Big(\frac{\partial P}{\partial\rho }\Big)_S }\Big|_{\rho=\rho_b, S=S_b}=
\sqrt{ \frac{\gamma P_b}{\rho_b} },
\end{equation}
the speed of sound. 

Note that 
$$
\mathcal{A}_b=(\mathfrak{a}_b|\bm{n}_b), \quad \mathcal{N}_b^2=
(\mathfrak{a}_b|\bm{n}_b)\Big(\frac{\nabla P_b}{\rho_b}\Big|\bm{n}_b\Big),
$$
where $\displaystyle \bm{n}_b=-\frac{\nabla \rho_b}{\|\nabla \rho_b\|}$, are the Schwarzschild discriminant, the square of the Brunt-V\"{a}is\"{a}l\"{a} 
frequency (local buoyancy frequency). \\

We note that the operator $\mathcal{L}_0$ can be written as 
\begin{align}
\mathcal{L}_0\bm{u}&=
\mathrm{grad}\Big[-\sigma_b\mathrm{div}(\rho_b\bm{u})+\sigma_b \rho_b(\bm{u}|\mathfrak{a}_b) \Big] + \nonumber \\
&+\sigma_b \Big[-\mathrm{div}(\rho_b\bm{u})\mathfrak{a}_b
+(\bm{u}|\mathfrak{a}_b)\nabla \rho_b \Big],
\end{align}
where
\begin{equation}
\sigma_b=\frac{\gamma P_b}{\rho_b^2}.
\end{equation}
Let us keep in mind that
\begin{align*}
&\frac{1}{C}\mathsf{d}^{\frac{1}{\gamma-1} }\leq \rho_b \leq {C}\mathsf{d}^{ \frac{1}{\gamma-1} }, \quad
\frac{1}{C}\mathsf{d}^{1+\frac{1}{\gamma-1} } \leq P_b \leq {C}\mathsf{d}^{ 1+\frac{1}{\gamma-1} },
\\
&
\frac{1}{C} \leq \frac{P_b}{\rho_b^{\gamma}} \leq C, 
\quad \frac{1}{C}\mathsf{d}^{ -\frac{2-\gamma}{\gamma-1} }\leq \sigma_b \leq {C}\mathsf{d}^{ -\frac{2-\gamma}{\gamma-1} }
\end{align*}
on $\mathfrak{R}_b$, where $\mathsf{d}=\mathsf{d}(\bm{x})=\mathrm{dist}(\bm{x}, \partial \mathfrak{R}_b)$\\

When $(\rho, S, \bm{v})=(\rho(t, \bm{x}), S(t,\bm{x}), \bm{v}(t, \bm{x}) )$ is a solution of the rotating Euler-Poisson equations \eqref{EPa}-\eqref{EPd} \eqref{EOS} near the stationary solution $(\rho_b, S_b, \bm{0})$, then the unknown variable $\bm{u}$ of the {\bf ELASO} means
\begin{equation}
\bm{u}(t,\bm{x})=\bm{\varphi}(t,\bm{x})-\bm{x} + \bm{u}^0(\bm{x}),
\end{equation}
where
$\bm{\varphi}$ is the flow of the velocity field $\bm{v}$ defined by
$$ \frac{\partial}{\partial t} \bm{\varphi}(t,\bm{x})=\bm{v}(t, \bm{\varphi}(t, \bm{x})),\quad \bm{\varphi}(0,\bm{x})=\bm{x},
$$ and $\bm{u}^0$ is supposed to satisfy
\begin{align*}
&\rho(0,\bm{x})-\rho_b(\bm{x})=-\mathrm{div}(\rho_b(\bm{x})\bm{u}^0(\bm{x})), \\
&S(0,\bm{x})-S_b(\bm{x})=-(\bm{u}^0(\bm{x})|\nabla S_b(\bm{x})).
\end{align*}\\

The formal integro-differential operator $\mathcal{L}$ considered on $C_0^{\infty}(\mathfrak{R}_b; \mathbb{C}^3)$ can be extended to a self-adjoint operator $\bm{L}$ in the Hilbert space
\begin{equation}
\mathfrak{H}=L^2(\mathfrak{R}_b, \rho_bd\bm{x}; \mathbb{C}^3)
\end{equation}
endowed with the norm
\begin{equation}
\|\bm{u}\|_{\mathfrak{H}}=\Big[ \int_{\mathfrak{R}_b}\|\bm{u}(\bm{x})\|\rho_b(\bm{x})d\bm{x} \Big]^{\frac{1}{2}}.
\end{equation}

Actually $\bm{L}$ is given by
\begin{align}
\mathsf{D}(\bm{L})&=\Big\{ \bm{u} \in \mathfrak{G}_0\quad\Big|\quad \mathcal{L}\bm{u} \in \mathfrak{H} \Big\}, \\
\bm{L}\bm{u}&=\mathcal{L}\bm{u}.
\end{align}
Here
\begin{equation}
\mathfrak{G}=\Big\{ \bm{u} \in \mathfrak{H}\quad\Big|\quad
\mathrm{div}(\rho_b\bm{u}) \in L^2(\mathfrak{R}_b, 
\sigma_bd\bm{x}; \mathbb{C})\quad \Big\}
\end{equation}
is a Hilbert space endowed with the norm 
\begin{equation}
\|\bm{u}\|_{\mathfrak{G}}=\Big[ \|\bm{u}\|_{\mathfrak{H}}^2+
\int_{\mathfrak{R}_b}|\mathrm{div}(\rho_b\bm{u})|^2\sigma_b(\bm{x})d\bm{x} \Big]^{\frac{1}{2}}
\end{equation}
and $\mathfrak{G}_0$ is the closure of $C_0^{\infty}(\mathfrak{R}_b; \mathbb{C}^3)$
in $ \mathfrak{G}$. \\

Details of mathematically rigorous discussion of the above described situation can be found in \cite{ELASO}. We use the following notations:

\begin{Notation}
Let $\mathsf{X}, \mathsf{Y}$ be Hilbert spaces. For an operator $T$ from a subspace of $\mathsf{X}$ into $\mathsf{Y}$, $\mathsf{D}(T)$ denotes the domain of $T$, 
\begin{align*}
&\mathsf{R}(T)=\mbox{the range of }T=\{ Tx |\quad x \in \mathsf{D}(T) \}, \\
&\mathsf{N}(T)=\mbox{the kernel of }T=\{ x \in \mathsf{D}(T) |\quad Tx =0_{\mathsf{Y}} \}.
\end{align*}
$\mathscr{B}(\mathsf{X};\mathsf{Y})$ denotes the Banach space of
all bounded linear operators from $\mathsf{X}$ into $\mathsf{Y}$: 
$$
|\|T\||_{\mathscr{B}(\mathsf{X};\mathsf{Y})}:=\sup\Big\{\|Tx\|_{\mathsf{Y}} \Big|\quad \|x\|_{\mathsf{X}}=1 \Big\}<\infty.
$$
$\mathscr{B}(\mathsf{X})=\mathscr{B}(\mathsf{X}; \mathsf{X})$

Let $T$ be an operator in $\mathsf{X}$ such that $\mathsf{D}(T)$ is dense in $\mathsf{X}$. 
\begin{align*}
&\mbox{\Pisymbol{psy}{82}}(T)=\mbox{the resolvent set of }T=\Big\{ \lambda \in \mathbb{C} \Big|\quad
\mathsf{N}(\lambda-T)=\{0_{\mathsf{X}}\}, 
(\lambda-T)^{-1} \in \mathscr{B}(\mathsf{X}) \Big\}, \\
&\mbox{\Pisymbol{psy}{83}}(T)=\mbox{the spectrum of }T=\mathbb{C} \setminus \mbox{\Pisymbol{psy}{82}}(T), \\
&\mbox{\Pisymbol{psy}{83}}_p(T)=\mbox{the set of all eigenvalues of }T=\Big\{ \lambda \in \mathbb{C} \Big|\quad 
\mathsf{N}(\lambda-T)\not=\{0_{\mathsf{X}}\} \Big\}.
\end{align*}
\end{Notation}

We are interested the structure of the spectrum $\mbox{\Pisymbol{psy}{83}}(\bm{L})$ of the self-adjoint operator $\bm{L}$. If $\mbox{\Pisymbol{psy}{83}}(\bm{L})=\mbox{\Pisymbol{psy}{83}}_p(\bm{L})$, where $\mbox{\Pisymbol{psy}{83}}_p(\bm{L})$ denotes the set of all eigenvalues of $\bm{L}$, then there is an orthonormal system of eigenfunctions, $(\bm{\phi}_n)_n$, $\bm{L} \bm{\phi}_n=\lambda_n\bm{\phi}_n$, which is complete in 
$\mathfrak{H}$. See, e.g., \cite[Theorem X.3.4]{DunfordS}. In this situation we have eigenfunction expansions
$\bm{u}=\sum c_n\bm{\phi}_n$ for $\forall \bm{u} \in \mathfrak{H}$, for which
$\bm{L}\bm{u}=\sum \lambda_nc_n\bm{\phi}_n$, and, if $\Omega=0$, the general solution of {\bf ELASO}
$$ \frac{\partial^2\bm{u}}{\partial t^2}+\bm{L}\bm{u}=\bm{0}$$
is given by
$$
\bm{u}(t, \bm{x})=\sum_n(c_n^+\bm{u}_n^+(t,\bm{x})+c_n^-\bm{u}_n^-(t,\bm{x})),
$$
where
$$
\bm{u}_n^{\pm}(t,\bm{x})=
\begin{cases}
e^{\pm\sqrt{\lambda_n}\mathrm{i}t}\bm{\phi}_n(\bm{x}) \quad (\lambda_n \geq 0) \\
e^{\pm\sqrt{|\lambda_n|}t}\bm{\phi}_n(\bm{x}) \quad (\lambda_n <0)
\end{cases}
$$
and
$$
c_n^{\pm}=
\begin{cases}
\frac{1}{2}
\Big((\bm{u}^0|\bm{\phi}_n)_{\mathfrak{H}}\pm\frac{1}{\sqrt{\lambda_n}\mathrm{i}}(\bm{v}^0|\bm{\phi}_n)_{\mathfrak{H}}\Big) \quad (\lambda_n >0) \\
\frac{1}{2}\Big((\bm{u}^0|\bm{\phi}_n)_{\mathfrak{H}}\pm\frac{1}{\sqrt{|\lambda_n|}}(\bm{v}^0|\bm{\phi}_n)_{\mathfrak{H}}\Big) \quad (\lambda_n <0) \\
\frac{1}{2}(\bm{u}^0|\bm{\phi}_n)_{\mathfrak{H}} \quad (\lambda_n=0),
\end{cases}
$$
with $\displaystyle \bm{u}^0(\bm{x})=\bm{u}(0,\bm{x}), \bm{v}^0(\bm{x})=\frac{\partial \bm{u}}{\partial t}(0,\bm{x})$.

However, if there are continuous spectra, that is, if $\mbox{\Pisymbol{psy}{83}}(\bm{L}) \setminus \mbox{\Pisymbol{psy}{83}}_p(\bm{L}) \not= \emptyset$, then the eigenfunction expansion does not work. In this sense, the question whether $\mbox{\Pisymbol{psy}{83}}(\bm{L}) = \mbox{\Pisymbol{psy}{83}}_p(\bm{L})$ or not is important. The aim of this study is concerned with this question, namely
\begin{Question}\label{Q.0}
Is it the case that $\mbox{\Pisymbol{psy}{83}}(\bm{L}) = \mbox{\Pisymbol{psy}{83}}_p(\bm{L})$ ?
\end{Question} 

Note that $\bm{L}$ is not of the Strum-Liouville type, say, with discrete spectrum, since 
the multiplicity of the eigenvalue $0$ is infinite, or, $\mathrm{dim}\mathsf{N}(\bm{L})=\infty$, and, therefore, the resolvent is not compact. In fact,

1) Suppose $\mathfrak{a}_b\not=0$. Then 
$$\bm{u}(\bm{x})=\frac{1}{\rho_b(\bm{x})}\mathfrak{a}_b(\bm{x})\times \nabla f(\bm{x}), $$
where $f \in C_0^{\infty}(\mathfrak{R}_b;\mathbb{R})$ is arbitrary, enjoys
$$\mathrm{div}(\rho_b\bm{u})=0,\quad (\bm{u}|\mathfrak{a}_b)=0, $$
since $\mathfrak{a}_b=-\frac{1}{\gamma \mathsf{C}_V}\nabla S_b$ satisfies
$\mathrm{rot}\mathfrak{a}_b=0$; Then $\bm{u} \in \mathsf{N}(\bm{L})$, since
$\delta \rho=0, \delta P=0$;

2) Suppose $\mathfrak{a}_b=0$. Then 
$$\bm{u}(\bm{x})=\frac{1}{\rho_b(\bm{x})}\mathrm{rot}\bm{f} (\bm{x}) $$
with arbitrary $\bm{f} \in C_0^{\infty}(\mathfrak{R}_b; \mathbb{R}^3)$ belongs to the kernel 
$\mathsf{N}(\bm{L})$.\\ 

 Here we note the following fact:\\

{\it 
Suppose that the background is isentropic, that is, $\mathfrak{a}_b=0$ everywhere.
Let $\bm{L}^G$ be the Friedrichs extension of the operator $\mathcal{L}\restriction C_0^{\infty}(\mathfrak{R}_b,; \mathbb{C}^3)$ in the functional space
$\mathfrak{G}=\{ \bm{u} \in \mathfrak{H} | \mathrm{div}(\rho_b\bm{u}) \in L^2(\mathfrak{R}_b, \sigma_bd\bm{x}; \mathbb{C})\}$. Then
 $ \mbox{\Pisymbol{psy}{83}}(\bm{L}^G)=
 \mbox{\Pisymbol{psy}{83}}_{\mathrm{p}}(\bm{L}^G)$, and
$ \mbox{\Pisymbol{psy}{83}}(\bm{L}^G)$ consists of $\{0\}$ and a sequence of eigenvalues $\lambda_n, n \in\mathbb{N},$ of finite multiplicities such that $\lambda_n\not=0, 
\lambda_n <\lambda_{n+1} \rightarrow  +\infty$ as $n \rightarrow \infty$.} \\

For proof see \cite[Sections III, IV ]{JJTM2020} for the case of spherically symmetric background for $\Omega=0$,
and \cite[Theorem 6]{ELASO}. Anyway we have $\mbox{\Pisymbol{psy}{83}}_p(\bm{L})=\mbox{\Pisymbol{psy}{83}}_p(\bm{L}^G)=\mbox{\Pisymbol{psy}{83}}(\bm{L}^G)$, and the Question is: Is 
$(\lambda -\bm{L})^{-1} \Big( \supset (\lambda- \bm{L}^G)^{-1} \Big) \in 
\mathscr{B}(\mathfrak{H})$
 when $\lambda \in \mbox{\Pisymbol{psy}{82}}(\bm{L}^G)$, that is, $ 
(\lambda-\bm{L}^G)^{-1} \in \mathscr{B}(\mathfrak{G})$ ?\\

On the other hand,\\

{\it Suppose $\mathfrak{a}_b\not=0$, $\Omega=0$. There can appear the so-called `$g$-mode'
$\{ \lambda_{-n} ; n \in \mathbb{N}\} \subset \mbox{\Pisymbol{psy}{83}}_p(\bm{L})$ such that $\lambda_{-n}>0$, $\lambda_{-n} \rightarrow 0$ as $ n \rightarrow \infty$. It is the case when
$\displaystyle \inf_{\mathfrak{R}_b}\frac{1}{r}\frac{dS_b}{dr} >0$,
or, $\displaystyle \inf_{\mathfrak{R}_b}\frac{\mathcal{N}_b^2}{r^2}>0$. }\\

For proof see \cite{TM2023}.\\

Sequential discussions are briefly as follows:

In Section 2 we introduce a first order system {\bf ELASO\ddag}, 
which is equivalent to the second order equation {\bf ELASO};

In Section 3 we introduce a first order system {\bf ELASO\ddag($\lozenge$)},
which is probably equivalent to the system {\bf ELASO\ddag}, and
discuss the esquivalence;

In Section 4 we analyze the generator $\bm{C}$ of the system
{\bf ELASO\ddag ($\lozenge$)} ;

In Section 5 we derive a sufficient condition for a complex number to be an essential spectrum of $\bm{C}$.

Thus, if the equivalnce between {\bf ELASO\ddag} and {\bf ELASO\ddag($\lozenge$)}
is justified exactly, then this is a sufficient condistion for a comlex number to
be an essential spectrum of the generator $\bm{A}$ of {\bf ELASO\ddag}.
This gives an answer to the 
{\bf Question 1}: Whether $\mbox{\Pisymbol{psy}{83}}(\bm{L}) = \mbox{\Pisymbol{psy}{83}}_p(\bm{L})$ or not.

\section{{\bf ELASO\ddag}}

The second order equation {\bf ELASO} \eqref{ELASO} is equivalent to the evolution equation, which we call {\bf ELASO\ddag}. :
\begin{equation}
\frac{\partial U}{\partial t}+\mathcal{A}U=0 \label{ELASO2}
\end{equation}
with
\begin{equation}
\mathcal{A}=
\begin{bmatrix}
\mathcal{B} & \mathcal{L} \\
\\
-I & O
\end{bmatrix}
\end{equation}
for the unknown
\begin{equation}
U=
\begin{bmatrix}
U^1 \\
U^2 \\
U^3 \\
U^4 \\
U^5 \\
U^6
\end{bmatrix}
=
\begin{bmatrix}
\bm{v} \\
\\
\bm{u}
\end{bmatrix}
=
\begin{bmatrix}
v^1 \\
v^2 \\
v^3 \\
u^1 \\
u^2 \\
u^3
\end{bmatrix},
\quad 
\bm{v}=
\begin{bmatrix}
v^1 \\
v^2 \\
v^3
\end{bmatrix}, \quad \bm{u}=
\begin{bmatrix}
u^1 \\
u^2 \\
u^3
\end{bmatrix}
.
\end{equation}\\

We consider the {\bf ELASO\ddag} in the Hilbert space
$\mathfrak{E}=\mathfrak{H}\times \mathfrak{G}_0$
with
the densely defined closed operator $\bm{A}$ in $\mathfrak{E}$,
$\mathsf{D}(\bm{A})=\mathfrak{G}_0 \times \mathsf{D}(\bm{L})$, $\bm{A}U=\mathcal{A}U$, namely,
$$
\bm{A}=\begin{bmatrix}
\bm{B} & \bm{L} \\
\\
-\bm{I} & O
\end{bmatrix},
$$
where
$\bm{B}: \bm{v} \mapsto \mathcal{B}\bm{v}$ is a bounded operator from $\mathfrak{H}$ onto $\mathfrak{H}$.

Note that $\mathsf{N}(\bm{A})=\{ \bm{0} \} \times \mathsf{N}(\bm{L})$.

Then, given $U_0 \in \mathsf{D}(\bm{A})$, the initial value problem
\begin{equation}
\frac{d U}{dt} +\bm{A}U=0, \quad U|_{t=0}=U_0
\end{equation}
admits a unique solution $U=U(t,\bm{x})$ in $C^1([0,+\infty[; \mathfrak{E})\cap C([0,+\infty[; \mathsf{D}(\bm{A}))$. And, for this $U(t,\bm{x})=(\bm{v}(t,\bm{x}), \bm{u}(t,\bm{x}))^{\top}$, the component $\bm{u}(t,\bm{x})$ is a solution of {\bf ELASO}\eqref{ELASO}
in $C^2([0,+\infty[; \mathfrak{H}) \cap C^1([0,+\infty[; \mathfrak{G}_0) \cap
C([0,+\infty[; \mathsf{D}(\bm{L}) )$, and $\bm{v}(t,\bm{x})=\partial \bm{u}(t,\bm{x})/\partial t$.

Note that 1) $\lambda \in \mbox{\Pisymbol{psy}{83}}(\bm{A})$ if and only if $\lambda^2-\lambda \bm{B}+\bm{L}$ does not have a bounded inverse, and 2) $ \lambda \in \mbox{\Pisymbol{psy}{83}}_p(\bm{A})$ if and only if there is $\bm{\phi} \in \mathsf{D}(\bm{L})$ such that $\bm{\phi} \not=0, 
(\lambda^2-\lambda \bm{B} + \bm{L})\bm{\phi}=\bm{0}$. Hence, when $\Omega=0, \bm{B}=O$, then it holds that 
$\mbox{\Pisymbol{psy}{83}}(\bm{L})=\mbox{\Pisymbol{psy}{83}}_p(\bm{L}) \Leftrightarrow 
\mbox{\Pisymbol{psy}{83}}(\bm{A})=\mbox{\Pisymbol{psy}{83}}_p(\bm{A})$,
since
$$\lambda \in \mbox{\Pisymbol{psy}{83}}(\bm{A})\ [\!( \in \mbox{\Pisymbol{psy}{83}}_p(\bm{A}) )\!] \quad
\Leftrightarrow \quad
-\lambda^2 \in \mbox{\Pisymbol{psy}{83}}(\bm{L})\ [\!( \in \mbox{\Pisymbol{psy}{83}}_p(\bm{L}) )\!] 
$$
for $\bm{B}=O$.\\

We consider

\begin{Question}\label{Q.1}
Is it the case that $\mbox{\Pisymbol{psy}{83}}(\bm{A}) = \mbox{\Pisymbol{psy}{83}}_p(\bm{A})$ ?
\end{Question}

When $\Omega=0, \bm{B}=O$, this {\bf Question \ref{Q.1} } is nothing but {\bf Question \ref{Q.0} }.

\section{ {\bf ELASO\ddag ($\lozenge$)} }

We transform {\bf ELASO\ddag}
to a first order system, which will be called {\bf ELASO\ddag ($\lozenge$)},
on the variables
\begin{equation}
W=
\begin{bmatrix}
W^1 \\
W^2 \\
W^3 \\
W^4 \\
W^5 \\
\end{bmatrix}
=
\begin{bmatrix}
\bm{v} \\
\\
\bm{w}
\end{bmatrix}=
\begin{bmatrix}
v^1 \\
v^2 \\
v^3 \\
w^1 \\
w^2
\end{bmatrix},
\quad
\bm{v}=
\begin{bmatrix}
v^1 \\
v^2 \\
v^3
\end{bmatrix},
\quad 
\bm{w}=
\begin{bmatrix}
w^1 \\
w^2\end{bmatrix}
,
\end{equation}
where

\begin{align}
\bm{w}=\mathcal{W}\bm{u}&=
\begin{bmatrix}
\mathcal{W}\restriction^1\bm{u} \\
\\
\mathcal{W}\restriction^2\bm{u}
\end{bmatrix} = \nonumber \\
&=
\begin{bmatrix}
\frac{\delta \rho}{\rho_b}-\frac{\delta P}{\gamma P_b} \\
\\
\frac{1}{{c_b} \rho_b}\delta P
\end{bmatrix}
=
\begin{bmatrix}
-(\bm{u}|\mathfrak{a}_b) \\
\\
-\frac{{c_b}}{\rho_b}\mathrm{div}(\rho_b\bm{u})+{c_b}(\bm{u}|\mathfrak{a}_b)
\end{bmatrix}.
\end{align}

The equation turns out to be
\begin{subequations}
\begin{align}
&\frac{\partial \bm{v}}{\partial t}+2\bm{\Omega}\times \bm{v} +\mathcal{L}^W\bm{w} =0 \\
&\frac{\partial w^1}{\partial t}+(\bm{v}|\mathfrak{a}_b)=0 \\
&\frac{\partial w^2}{\partial t}+\frac{{c_b}}{\rho_b}\mathrm{div}(\rho_b \bm{v})-{c_b} (\bm{v}|\mathfrak{a}_b)=0
\end{align}
\end{subequations}
with
\begin{align}
\mathcal{L}^W\bm{w}&=\frac{1}{\rho_b}\nabla({c_b} \rho_b w^2)-\frac{\nabla P_b}{{c_b} \rho_b}({c_b} w^1+w^2) 
-4\pi\mathsf{G}\nabla \mathcal{K}\Big[\rho_b w^1+\frac{\rho_b}{{c_b}}w^2\Big] \nonumber \\
&=\mathcal{L}^W_{01}w^1+ \mathcal{L}^W_{02}w^2
+4\pi\mathsf{G}\mathcal{L}_1^W\bm{w}, \\
&\mathcal{L}^W_{01}w^1=-{c_b}^2\frac{\nabla k_1}{k_1} w^1 \\
&\mathcal{L}^W_{02}w^2=
{c_b}\frac{1}{k_2}\nabla (k_2 w^2) \\
&\mathcal{L}_1^W\bm{w}=-
\nabla \mathcal{K}\Big[\rho_b w^1+\frac{\rho_b}{{c_b}}w^2\Big].
\end{align}
Here we have introduced
the coefficients
\begin{align}
k_1&=P_b^{\frac{1}{\gamma}}=\rho_b \cdot E, \nonumber \\
k_2&=\sqrt{\gamma}\rho_b^{\frac{1}{2}}P_b^{-\frac{2-\gamma}{2\gamma}}
={c_b} \cdot \frac{1}{E}, \nonumber \\
k_3&=\rho_bP_b^{-\frac{1}{\gamma}}=
\frac{1}{E},
\end{align}
where 
\begin{equation}
E=e^{S_b/\gamma \mathsf{C}_V}.
\end{equation}
Note that 
both $E$ and $\displaystyle \frac{1}{E}$ belong to $C^{3,\alpha}(\mathfrak{R}_b \cup \partial \mathfrak{R}_b)$
and
\begin{equation}
\frac{\nabla E}{E}= -\mathfrak{a}_b.
\end{equation}\\


Note that
\begin{equation}
k_1k_2={c_b} \rho_b,\quad k_2={c_b} k_3,
\quad \frac{\nabla k_3}{k_3}=\mathfrak{a}_b, 
\quad \frac{{c_b} k_2}{k_1}=\sigma_b \frac{1}{E^2}.
\end{equation}

As for the behavior at the vacuum boundary of the coefficients, we note 
\begin{align}
&0<\frac{1}{C} \mathsf{d}^{\frac{1}{\gamma-1}} \leq k_1 
\leq {C} \mathsf{d}^{\frac{1}{\gamma-1}} ,
\quad
0<\frac{1}{C}\mathsf{d}^{\frac{\gamma-1}{2}} \leq k_2 \leq {C}\mathsf{d}^{\frac{\gamma-1}{2}}, \nonumber \\
&
0<\frac{1}{C}\leq k_3 \leq C
\end{align}
on $\mathfrak{R}_b$, where $\mathsf{d}=\mathrm{dist}(\cdot, \partial \mathfrak{R}_b)$,
and
\begin{equation}
\frac{{c_b}^2}{ k_1}\nabla k_1 \Big(=-\frac{\nabla P_b}{\rho_b} \Big) \quad\mbox{and}\quad \frac{1}{k_3}\nabla k_3 \Big(=\mathfrak{a}_b \Big) \quad
\in C^{0,\alpha}(\mathfrak{R}_b\cup \partial \mathfrak{R}_b ).
\end{equation}\\

We shall often use the relation
\begin{align}
\frac{c_b}{k_1}\mathrm{div}(k_1\bm{v})&=
\frac{1}{\sqrt{\rho_b}}\cdot\sqrt{\sigma_b}\cdot \frac{1}{E}\mathrm{div}(E\cdot \rho_b
\bm{v}) \nonumber \\
&=\frac{c_b}{\rho_b}\mathrm{div}(\rho_b\bm{v})-c_b
(\bm{v}|\mathfrak{a}_b).
\end{align} 
Therefore we can write
\begin{equation}
\mathcal{W}\bm{u}=
\begin{bmatrix}
-\frac{1}{k_3}(\bm{u}| \nabla k_3) \\
\\
-\frac{c_b}{k_1}\mathrm{div}(k_1\bm{u})
\end{bmatrix},
\end{equation}
and
\begin{align}
\mathcal{L}\bm{u}&=\mathcal{L}^W\mathcal{W}\bm{u} = \nonumber \\
&=(\bm{u}|\mathfrak{a}_b) \frac{\nabla P_b}{\rho_b}+
E\mathrm{grad}\Big[-\frac{\sigma_b}{E}\mathrm{div}\Big(E\rho_b\bm{u}\Big)\Big] 
+4\pi\mathsf{G}\nabla \mathcal{K}[\mathrm{div}(\rho_b\bm{u})],
\end{align}
where we note
$\displaystyle \frac{\nabla P_b}{\rho_b}
\in C^{1,\alpha}(\mathfrak{R}_b\cup \partial\mathfrak{R}_b; \mathbb{R}^3)$.\\

The system to be considered is
\begin{equation}
\frac{\partial W}{\partial t}+\mathcal{C} W=0, \label{ELASO2loz}
\end{equation}
where
\begin{equation}
\mathcal{C}=
\begin{bmatrix}
\mathcal{B} & \mathcal{L}^W \\
\\
-\mathcal{W} & 0
\end{bmatrix} 
.
\end{equation}

While $\mathcal{C}, \mathcal{L}^W, \mathcal{W}$ are formal integro-differential operators, we are going to fix the idea on operators in the space
\begin{equation}
\mathfrak{E}^W= \mathfrak{H} \times \mathfrak{h}^2.
\end{equation}
Here and hereafter we denote
$\mathfrak{h}^D=L^2(\mathfrak{R}_b, \rho_b d\bm{x}; \mathbb{C}^D)$ for 
$D=1,2,3,4,5$,
while $\mathfrak{H}=\mathfrak{h}^3$.\\

First $\bm{B}: \bm{v} \mapsto \mathcal{B}\bm{v}=2{\Omega}
\begin{bmatrix}
-v^2 \\
v^1 \\
0
\end{bmatrix}
$ is an operator in $\mathscr{B}(\mathfrak{H})$.\\

Next $\bm{W}: \bm{u} \mapsto \mathcal{W}\bm{u}$ is an operator in $\mathscr{B}(\mathfrak{G}_0 ; \mathfrak{h}^2)$.\\

As for $\mathcal{L}^W$, we consider $\bm{L}^W$ defined by
\begin{equation}
\mathsf{D}(\bm{L}^W)=\mathfrak{h}^1\times \mathfrak{f}, \qquad
\bm{L}^W\bm{w}=\mathcal{L}^W\bm{w}.
\end{equation}
Here
\begin{equation}
\mathfrak{f}=\Big\{ 
w \in \mathfrak{h}^1 \Big|\quad \mathcal{L}_{02}^Ww=
\frac{c_b}{k_2}\mathrm{grad}(k_2 w) \in \mathfrak{H} \Big\}.
\end{equation}

Since $\mathfrak{f}$ is dense in $\mathfrak{h}^1$, $\mathsf{D}(\bm{L}^W)$ is dense in $\mathfrak{h}^2$.\\

We claim 

\begin{Lemma}
The operator $\bm{L}^W$ densely defined in $\mathfrak{h}^2$ into $\mathfrak{H}$ is a closed operator.
\end{Lemma}

Proof. Let us consider a sequence $(\bm{w}_n)_n$ in $\mathsf{D}(\bm{L}^W)$ such that $\bm{w}_n \rightarrow \bm{w}$ in $\mathfrak{h}^1$ and $\bm{L}^W\bm{w}_n \rightarrow \bm{f}$ in $\mathfrak{H}$. We want to deduce $w^2 \in \mathfrak{f}$. Look at
$$\mathcal{L}^W\bm{w}_n=\mathcal{L}_{01}^Ww_n^1+ \mathcal{L}_{02}^Ww_n^2+4\pi\mathsf{G}\mathcal{L}_1^W\bm{w}_n.
$$
Since $\mathcal{L}^W\bm{w}_n, \mathcal{L}_{01}^Ww_n^1, 4\pi\mathsf{G}\mathcal{L}_1^W\bm{w}_n$ converge to
$\bm{f}, \mathcal{L}_{01}^Ww^1, 4\pi\mathsf{G}\mathcal{L}_1^W\bm{w}$ as $n \rightarrow \infty$, we see
$\mathcal{L}_{02}^Ww_n^2$ converges to $\bm{f}_{02}$, where $\bm{f}_{02}:=\bm{f} - \mathcal{L}_{01}^Ww^1 - 4\pi\mathsf{G}\mathcal{L}_1^W\bm{w}$. Then any test function $\bm{\varphi} \in C_0^{\infty}(\mathfrak{R}_b; \mathbb{C}^3)$ enjoys
\begin{align*}
&-\int_{\mathfrak{R}_b}k_2w^2\mathrm{div}\Big(\frac{c_b}{k_2}\bm{\varphi}\Big)^*d\bm{x} =
\lim_n\Big[
-\int_{\mathfrak{R}_b}k_2w_n^2\mathrm{div}\Big(\frac{c_b}{k_2}\bm{\varphi}\Big)^*d\bm{x} \Big] \\
&=\lim_n\int_{\mathfrak{R}_b}\frac{c_b}{k_2}\Big(\mathrm{grad}(k_2 w_n^2)\Big| \bm{\varphi}\Big)d\bm{x} =
\lim_n\int_{\mathfrak{R}_b}(\mathcal{L}_{02}^Ww_n^2\Big|\bm{\varphi})d\bm{x} \\
&=\int_{\mathfrak{R}_b}(\bm{f}_{02}|\bm{\varphi})d\bm{x},
\end{align*}
since $\mathcal{L}_{02}^Ww_n^2 \rightarrow \bm{f}_{02} $ in $L^2(\mathrm{supp}[\bm{\varphi}])$ for
$\sup_{\mathrm{supp}[\bm{\varphi}]}\frac{1}{\rho_b} <\infty$.
Therefore $\mathcal{L}_{02}^Ww^2=\bm{f}_{02}$ in the distribution sense and $w^2 \in \mathfrak{f}$; Hence
$\bm{w} \in \mathfrak{h}^1\times \mathfrak{f}, \bm{L}^W\bm{w}=\bm{f}$. $\square$\\

Note that $\bm{W}\mathsf{D}(\bm{L}) \subset\mathsf{D}(\bm{L}^W)$ and
$$ \bm{L}^W\bm{W}\bm{u}=\bm{L}\bm{u} \quad\mbox{for}\quad \bm{u} \in \mathsf{D}(\bm{L}).
$$

\begin{Remark}
We cannot claim that $\bm{W}\mathsf{D}(\bm{L}) $, or $\bm{W}\mathfrak{G}_0$, is dense in $\mathfrak{h}^2$. In fact, when 
$\mathfrak{a}_b=0$, $\bm{W}\mathfrak{G}_0 \subset \{0\} \times \mathfrak{h}^1$, which is not dence in
$\mathfrak{h}^2$. Even when $\mathfrak{a}_b \not=0$, we see
$\bm{W}\mathfrak{G}_0 \subset \mathfrak{h}^1 \times \mathfrak{h}_C^1$, 
where
$$
\mathfrak{h}_C^1=\Big\{ w \in \mathfrak{h}^1 \Big|\quad
\int_{\mathfrak{R}_b}\frac{k_1}{c_b}wd\bm{x}=\Big(\frac{1}{k_2}\Big| w\Big)_{\mathfrak{h}^1}=0 \Big\},
$$
which is a closed subspace of $\mathfrak{h}^1$ with codimension $1$
so that it is not dense in $\mathfrak{h}^1$.
We have not yet found a neat characterization of $\bm{W}\mathfrak{G}_0$, or of $\bm{W}\mathsf{D}(\bm{L}) $, as a subspace of $\mathfrak{h}^2$.
\end{Remark}

We claim
\begin{Lemma}
It holds that
\begin{equation}
\Big\{ \bm{u} \in \mathfrak{G}_0 \Big|\quad \bm{W}\bm{u} \in \mathsf{D}(\bm{L}^W) \Big\}
\subset
\mathsf{D}(\bm{L}).
\end{equation}
\end{Lemma}

Proof. Let $\bm{u} \in \mathfrak{G}_0$ and $\bm{W}\bm{u} \in \mathsf{D}(\bm{L}^W)=\mathfrak{h}^1\times \mathfrak{f}$. We want to deduce $\mathcal{L}\bm{u} \in \mathfrak{H}$. Since
$$\mathcal{L}\bm{u}=\mathcal{L}_{01}^Ww^1+\mathcal{L}_{02}^Ww^2+4\pi\mathsf{G}\mathcal{L}_1^W\bm{w}
$$
for $\bm{w}=\mathcal{W}\bm{u}$, it is sufficient to deduce
$\mathcal{L}_{02}^W w^2 \in \mathfrak{H}$. But it is the case since $w^2 \in \mathfrak{f}$. 
Hence $\bm{u} \in \mathsf{D}(\bm{L})$. $\square$\\

We fix our idea by putting
\begin{equation}
\bm{C} =
\begin{bmatrix}
\bm{B} & \bm{L}^W \\
\\
-\bm{W} & O
\end{bmatrix},
\quad 
\mathsf{D}(\bm{C})=\mathfrak{G}_0\times \mathsf{D}(\bm{L}^W).
\end{equation}

The domain $\mathsf{D}(\bm{C})$ is dense in $\mathfrak{E}^W=\mathfrak{H} \times \mathfrak{h}^2$, since $\mathsf{D}(\bm{L}^W)$ is dense in $\mathfrak{h}^2$.\\

Since $\bm{L}^W$ is closed, we can claim that the operator $\bm{C}$ is 
a densely defined closed operator in $\mathfrak{E}^W$.\\

We put
$$\widetilde{\bm{W}}=
\begin{bmatrix}
I & O \\
\\
O & \bm{W}
\end{bmatrix}
: 
\mathfrak{E} \rightarrow \mathfrak{E}^W :
\begin{bmatrix}
\bm{v} \\
\\
\bm{u}
\end{bmatrix}
\mapsto
\begin{bmatrix}
\bm{v} \\
\\
\bm{W}\bm{u}
\end{bmatrix}
.
$$

\begin{Lemma}
If $U=U(t,\bm{x})
=\begin{bmatrix}
\bm{u}(t,\bm{x}) \\
\\
\bm{v}(t,\bm{x})
\end{bmatrix}$ is a solution of {\bf ELASO\ddag} in
$C^1([0,+\infty[; \mathfrak{E})\cap C([0,+\infty[; \mathsf{D}(\bm{A}))$, then the corresponding
$$W=W(t,\bm{x}):=\widetilde{\bm{W}}
U(t,\bm{x}) 
=\begin{bmatrix}
\bm{v}(t,\bm{x}) \\
\\
\bm{W}\bm{u}(t,\bm{x})
\end{bmatrix}
$$ turns out to be a solution of {\bf ELASO\ddag($\lozenge$)} in
$C^1([0,+\infty[; \mathfrak{E}^W )
\cap C([0,+\infty[; \mathsf{D}(\bm{C}))$. 
\end{Lemma}

Let us note 

\begin{Lemma}
It holds 
that $\mbox{\Pisymbol{psy}{83}}_p(\bm{C}) \setminus \{0\} =\mbox{\Pisymbol{psy}{83}}_p(\bm{A})\setminus \{0\}$.
\end{Lemma}

Proof. Let $\lambda \in \mbox{\Pisymbol{psy}{83}}_p(\bm{A}), \lambda \not=0$. Then there exists
$
\begin{bmatrix}
\bm{v}_0 \\
\\
\bm{u}_0
\end{bmatrix}
\in \mathfrak{G}_0 \times \mathsf{D}(\bm{L}), \not=
\begin{bmatrix}\bm{0} \\
\\
\bm{0}
\end{bmatrix}
$
such that
$$\bm{B}\bm{v}_0+\bm{L}\bm{u}_0 -\lambda \bm{v}_0=\bm{0},\quad
-\bm{v}_0-\lambda \bm{u}_0 =\bm{0}.
$$
Put $\bm{w}_0=\bm{W}\bm{u}_0 \in \mathsf{D}(\bm{L}^W)$. Then
$\bm{L}^W\bm{w}_0=\bm{L}\bm{u}_0$ and
$$
\bm{B}\bm{v}_0+\bm{L}^W\bm{w}_0 -\lambda \bm{v}_0 =\bm{0},\quad -\bm{W}\bm{v}_0 -\lambda \bm{w}_0=\bm{0}.
$$That is,
$$(\bm{C}-\lambda)
\begin{bmatrix}
\bm{v}_0 \\
\\
\bm{w}_0
\end{bmatrix}
=\begin{bmatrix}
\bm{0} \\
\\
\bm{0}
\end{bmatrix}
, $$
and $\bm{v}_0 \not=\bm{0}$, since, otherwise $\bm{u}_0=-\frac{1}{\lambda}\bm{v}_0 =\bm{0}$, contradicting to 
$
\begin{bmatrix}
\bm{v}_0 \\
\bm{u}_0
\end{bmatrix}
\not=
\begin{bmatrix}
\bm{0}\\
\bm{0}
\end{bmatrix}
$. Hence 
$
\begin{bmatrix}
\bm{v}_0 \\
\bm{w}_0
\end{bmatrix}
\not=
\begin{bmatrix}
\bm{0}\\
\bm{0}
\end{bmatrix}
$
and
$\lambda \in \mbox{\Pisymbol{psy}{83}}_p(\bm{C})$.

Let $\lambda \in \mbox{\Pisymbol{psy}{83}}_p(\bm{C}), \lambda \not=0 $.
Then there exists 
$\begin{bmatrix}
\bm{v}_0 \\
\\
\bm{w}_0
\end{bmatrix}
\in \mathfrak{G}_0\times \mathsf{D}(\bm{L}^W), \not=
\begin{bmatrix}
\bm{0}\\
\\
\bm{0}
\end{bmatrix}
$
such that
$$
\bm{B}\bm{v}_0+\bm{L}^W\bm{w}_0 -\lambda \bm{v}_0 =\bm{0},
\quad
-\bm{W}\bm{v}_0 -\lambda\bm{w}_0=\bm{0}.
$$
Since $\bm{w}_0 \in \mathsf{D}(\bm{L}^W)$, we have
$\bm{W}\bm{v}_0=-\lambda \bm{w}_0 \in \mathsf{D}(\bm{L}^W)$, therefore,
$\bm{v}_0$ being in $\mathfrak{G}_0$, $\bm{v}_0 \in \mathsf{D}(\bm{L})$ and
$$\bm{L}^W\bm{w}_0=-\frac{1}{\lambda}\bm{L}^W\bm{W}\bm{w}_0 =-\frac{1}{\lambda}
\bm{L}\bm{v}_0. $$
Then we have
$$
(\bm{A}-\lambda)
\begin{bmatrix}
\bm{v}_0 \\
\\
-\frac{1}{\lambda}\bm{v}_0
\end{bmatrix}
=
\begin{bmatrix}
\bm{0}\\
\\
\bm{0}
\end{bmatrix}
.
$$
Moreover $\bm{v}_0 \not=\bm{0}$, since, otherwise $\bm{w}_0=-\frac{1}{\lambda}\bm{W}\bm{v}_0=\bm{0}$,
contradicting to
$\begin{bmatrix}
\bm{v}_0 \\
\bm{w}_0
\end{bmatrix}
\not=
\begin{bmatrix}
\bm{0}\\
\bm{0}
\end{bmatrix}
$.
Therefore 
$\begin{bmatrix}
\bm{v}_0 \\
-\frac{1}{\lambda}\bm{v}_0
\end{bmatrix}
\not=
\begin{bmatrix}
\bm{0}\\
\bm{0}
\end{bmatrix}
$
and $\lambda \in \mbox{\Pisymbol{psy}{83}}_p(\bm{A})$. $\square$\\

We consider

\begin{Supposition}\label{Q.X}
It holds that $\mbox{\Pisymbol{psy}{83}}(\bm{C}) \subset \mbox{\Pisymbol{psy}{83}}(\bm{A})$, namely, $\mbox{\Pisymbol{psy}{82}}(\bm{A}) \subset \mbox{\Pisymbol{psy}{82}}(\bm{C})$.
\end{Supposition}

Later we shall show that, when $\mathfrak{a}_b \not=\bm{0}$,
there can exist $\alpha_{\pm}, \alpha_- <0 < \alpha_+,$ such that
$$ \Big\{ \lambda \Big|\quad -\lambda^2 \in [\alpha_-, \alpha_+] \Big\} 
\subset \mbox{\Pisymbol{psy}{83}}(\bm{C}).
$$
Therefore, if {\bf Supposition \ref{Q.X}} is the case, then 
$$ \Big\{ \lambda \Big|\quad -\lambda^2 \in [\alpha_-, \alpha_+] \Big\} 
\subset \mbox{\Pisymbol{psy}{83}}(\bm{A}),
$$
and $\mbox{\Pisymbol{psy}{83}}(\bm{A}) \setminus \mbox{\Pisymbol{psy}{83}}_p(\bm{A})\not=\emptyset $, that is,
{\bf Question \ref{Q.1}} is negatively answered, and {\bf Question \ref{Q.0}} is negatively answered when $\Omega=0$.\\

Let us observe what is the point in view of {\bf Supposition \ref{Q.X}}. Let $\lambda 
\in \mbox{\Pisymbol{psy}{82}}(\bm{A})$. Then $\lambda\not=0$ and
$(\lambda^2-\lambda \bm{B}+\bm{L})^{-1} \in \mathscr{B}(\mathfrak{H})$. If we want to show that $\lambda \in \mbox{\Pisymbol{psy}{82}}(\bm{C})$, we have to find
$\bm{v} \in \mathfrak{G}_0, \bm{w} \in \mathsf{D}(\bm{L}^W)$ such that
$$
\Big(\bm{C}-\lambda\Big)
\begin{bmatrix}
\bm{v} \\
\\
\bm{w}
\end{bmatrix}
=
\begin{bmatrix}
\bm{f} \\
\\
\bm{g}
\end{bmatrix}
,
$$
or
\begin{align*}
\bm{B}\bm{v}+\bm{L}^W\bm{w} -\lambda\bm{v}&=\bm{f} \\
-\bm{W}\bm{v}-\lambda\bm{w}&=\bm{g}
\end{align*}
for given $\bm{f} \in \mathfrak{H}, \bm{g} \in \mathfrak{h}^2$ with
$$\|\bm{v}\|_{\mathfrak{H}}^2+\|\bm{w}\|_{\mathfrak{h}^2}^2
\leq C^2 \Big[\|\bm{f}\|_{\mathfrak{H}}^2+\|\bm{g}\|_{\mathfrak{h}^2}^2\Big].
$$

First we claim the existence of $(\bm{C}-\lambda)^{-1}$. 

(Proof. Let
$$ (\bm{C}-\lambda)
\begin{bmatrix}
\bm{v} \\
\bm{w}
\end{bmatrix}
=\begin{bmatrix}
\bm{0} \\
\bm{0}
\end{bmatrix}
, \quad
\bm{v} \in \mathfrak{G}_0, \bm{w} \in \mathsf{D}(\bm{L}^W).
$$
Then 
$$\bm{B}\bm{v}+\bm{L}^W\bm{w}-\lambda \bm{v} =\bm{0}, \quad
-\bm{W}\bm{v}-\lambda \bm{w}=\bm{0}.
$$
Since $\lambda\not=0$, $\bm{w}=-\frac{1}{\lambda}\bm{W}\bm{v} \in \mathsf{D}(\bm{L}^W)$. Since $\bm{v} \in \mathfrak{G}_0$, we see $\bm{v} \in \mathsf{D}(\bm{L})$ and
$$\bm{L}^W\bm{w}=-\frac{1}{\lambda}\bm{L}\bm{v}.
$$
Then we have
$$-\lambda \bm{B}\bm{v}+\bm{L}\bm{v}+\lambda^2\bm{v}=\bm{0}, \quad \bm{v} \in \mathsf{D}(\bm{L}).
$$
Therefore $\bm{v}=\bm{0}$, and $\bm{w} =\bm{0}$. $\square$ )

Next we claim 

\begin{Lemma}
It holds that $\mathfrak{H} \times \mathsf{D}(\bm{L}^W) \subset
\mathsf{R}(\bm{C}-\lambda)$ and
$$(\bm{C}-\lambda)^{-1}\restriction \mathfrak{H}\times \mathsf{D}(\bm{L}^W) \in \mathscr{B}( \mathfrak{H}\times \mathsf{D}(\bm{L}^W)
,\mathfrak{E}^W).$$
\end{Lemma}

Proof. We want to find $\bm{v} \in \mathfrak{G}_0, \bm{w} \in \mathsf{D}(\bm{L}^W)$ such that
$$\bm{B}\bm{v}+\bm{L}^W\bm{w}-\lambda\bm{v}=\bm{f},\quad
-\bm{W}\bm{v}-\lambda\bm{w} =\bm{g}$$
for given $\bm{f} \in \mathfrak{H}, \bm{g} \in \mathsf{D}(
\bm{L}^W)$. But it is possible by solving
\begin{align*}
\bm{v}&=-(\lambda^2-\lambda\bm{B}+\bm{L})^{-1}(\lambda\bm{f}+\bm{L}^W\bm{g}) \\
\bm{w}&=\frac{1}{\lambda}\Big[
\bm{W}
(\lambda^2-\lambda\bm{B}+\bm{L})^{-1}(\lambda\bm{f}+\bm{L}^W\bm{g})
-\bm{g}\Big],
\end{align*}
since we are supposing $\bm{g} \in \mathsf{D}(\bm{L}^W)$. The norm
$\|\bm{v}\|_{\mathfrak{H}}, \|\bm{w}\|_{\mathfrak{h}^2}$ can be bounded by
$|\lambda|\|\bm{f}\|_{\mathfrak{H}} +\|\bm{g}\|_{\mathsf{D}(\bm{L}^W) } $,
since $(\lambda^2-\lambda\bm{B}+\bm{L})^{-1} \in \mathscr{B}(\mathfrak{H}, \mathfrak{G}_0)$
and $\bm{W} \in \mathscr{B}(\mathfrak{G}_0, \mathfrak{h}^2)$. $\square$ \\

Consequently $\mathsf{R}(\bm{C}-\lambda)$ is dense in $\mathfrak{E}^W$ and
the validity of 
the {\bf Supposition} reduces to the boundedness of $(\bm{C}-\lambda)^{-1}$
with respect to the norm $\|\cdot\|_{\mathfrak{h}^5}$. In other words, we are fronted with the alternative either $\lambda \in \mbox{\Pisymbol{psy}{82}}(\bm{C})$ or $\lambda \in \mbox{\Pisymbol{psy}{83}}_c(\bm{C})$ ( the continuous spectrum of $\bm{C}$ ),
for $\lambda \in \mbox{\Pisymbol{psy}{82}}(\bm{A})$ given. . We do not know whether the lattar possibility is excludable or not.


\section{Analysis of the operator of {\bf ELASO\ddag ($\lozenge$)} }

In order to analyze the operator $\mathcal{C}$, we decompose it as

\begin{equation}
\mathcal{C}=
\begin{bmatrix}
\mathcal{B} & \mathcal{L}^W \\
\\
-\mathcal{W} & 0
\end{bmatrix} 
=
{\mathcal{F}}+{\mathcal{H}}+{\mathcal{G}},
\end{equation}
where
\begin{align}
{\mathcal{F}}&=
\begin{bmatrix}
O & \bm{0} & \mathcal{L}_{02}^W \\
\bm{0}^{\top} & 0 & 0 \\
-\mathcal{W}\restriction^2 & 0 & 0
\end{bmatrix}
\nonumber \\
&=
\begin{bmatrix}
0 & 0 & 0 & 0 & \frac{{c_b}}{k_2}\partial_1(k_2 \cdot) \\
\\
0 & 0 & 0 & 0 & \frac{{c_b}}{k_2}\partial_2(k_2 \cdot) \\
\\
0 & 0 & 0 & 0 & \frac{{c_b}}{k_2}\partial_3(k_2 \cdot) \\
\\
0 & 0 & 0 & 0 & 0 \\
\\
\frac{{c_b}}{k_1}\partial_1(k_1\cdot) & \frac{{c_b}}{k_1}\partial_2(k_1\cdot) & \frac{{c_b}}{k_1}\partial_3(k_1\cdot)& 0 & 0
\end{bmatrix}, 
\end{align}
\begin{align}
{\mathcal{H}}&=
\begin{bmatrix}
\mathcal{B} & \mathcal{L}_{01}^W & \bm{0} \\
-\mathcal{W}\restriction^1 & 0 & 0 \\
\bm{0}^{\top} & 0 & 0 
\end{bmatrix}
\nonumber \\
&=
\begin{bmatrix}
0 & -2\Omega & 0 & -\frac{{c_b}^2}{k_1}\partial_1k_1 & 0 \\
\\
2\Omega & 0 & 0 & -\frac{{c_b}^2}{k_1}\partial_2k_1 & 0 \\
\\
0 & 0 & 0 & -\frac{{c_b}^2}{k_1}\partial_3k_1 & 0 \\
\\
\frac{1}{k_3}\partial_1 k_3 & \frac{1}{k_3}\partial_2 k_3 & \frac{1}{k_3}\partial_3 k_3 & 0 & 0 \\
\\
0 & 0 & 0 & 0 & 0
\end{bmatrix}, 
\end{align}
\begin{align}
{\mathcal{G}}&=
-4\pi\mathsf{G}
\begin{bmatrix}
O & \nabla\mathcal{K}[\rho_b\cdot] & \nabla\mathcal{K}[\frac{\rho_b}{c_b}\cdot] \\
\bm{0}^{\top} & 0 & 0 \\
\bm{0}^{\top} & 0 & 0
\end{bmatrix}
\nonumber \\
&=
-4\pi\mathsf{G}
\begin{bmatrix}
0 & 0 & 0 & \partial_1\mathcal{K}[\rho_b\cdot] & \partial_1 \mathcal{K}\Big[\frac{\rho_b}{{c_b}}\cdot\Big] \\
\\
0 & 0 & 0 & \partial_2\mathcal{K}[\rho_b\cdot] & \partial_2 \mathcal{K}\Big[\frac{\rho_b}{{c_b}}\cdot\Big] \\
\\
0 & 0 & 0 & \partial_3\mathcal{K}[\rho_b\cdot] & \partial_3 \mathcal{K}\Big[\frac{\rho_b}{{c_b}}\cdot\Big] \\
\\
0 & 0 & 0 & 0 & 0 \\
\\
0 & 0 & 0 & 0 & 0
\end{bmatrix}.
\end{align}

Here $\partial_j$ stands for $\displaystyle \frac{\partial}{\partial x^j}, j=1,2,3$.\\

First we claim 

\begin{Lemma}
The operator $\bm{G}$ defined as $ \mathsf{D}(\bm{G})=\mathfrak{E}^W, \bm{G}W=\mathcal{G}W$ is a compact operator. 
\end{Lemma}

Proof.
We see that $\bm{w} \mapsto g=\rho_b w^1+\frac{\rho_b}{{c_b}}w^2$ is a continuous mapping from $\mathfrak{h}^2$ into $L^2(\mathfrak{R}_b, \frac{\gamma P_b}{\rho_b^2}d\bm{x})$, which is continuously imbedded into $L^2(d\bm{x})$. On the other hand, $g \mapsto \mathcal{K}[g]$ is continuous from $L^2(d\bm{x})$ into $ H^2(d\bm{x})$, which is continuously imbedded into $H^1(d\bm{x})$. (
\cite[p.230, Theorem 9.9]{GilbergT}. ) Hence $g \mapsto \mathrm{grad}\mathcal{K}[g]$ is a compact operator from $L^2(d\bm{x})$ into $L^2(d\bm{x})$, which is continuously imbedded into $\mathfrak{h}^5$. Hence $\bm{G}$ is a compact operator in 
$\mathfrak{E}^W$. $\square$. \\

Next, 
$\mathcal{H}$ is a multiplication operator and its coefficients,
$\displaystyle 2\Omega, \frac{{c_b}^2}{k_1}\nabla k_1, \frac{1}{k_3}\nabla k_3$
all belong to $C^{0,\alpha}(\mathfrak{R}\cup\partial\mathfrak{R}) \subset L^{\infty}(\mathfrak{R})$. Therefore the operator $\bm{H}$ defined as
$\mathsf{D}(\bm{H})=\mathfrak{E}^W, \bm{H}W=\mathcal{H}W$ is a bounded operator in $\mathfrak{E}^W$.\\

\begin{Remark}
For our case $-\frac{{c_b}^2}{k_1}\mathrm{grad} k_1$ is bounded near the vacuum boundary, but $-\frac{{c_b}}{k_1}\mathrm{grad}k_1$ is not. This is the reason why we use the variables
$w_1=\frac{\delta\rho_b}{\rho_b}-\frac{\delta P}{\gamma P_b}, w^2=\frac{\delta P}{{c_b} \rho_b}$ instead of
the variables 
$m={c_b}\Big(\frac{\delta\rho_b}{\rho_b}-\frac{\delta P}{\gamma P_b}\Big), n=\frac{\delta P}{{c_b} \rho_b}$, so called Eckart variables, used in \cite{LL1992}, \cite{FLMM}.
\end{Remark}

Next we look at the operator $\bm{F}$:

$$ \mathsf{D}(\bm{F})=\mathsf{D}(\bm{C})=\mathfrak{G}_0 \times (\mathfrak{h}^1\times \mathfrak{f}),\quad \bm{F}W=\mathcal{F}W,$$
which is a densely defined closed operator in $\mathfrak{E}^W$.

We are considering
\begin{equation}
\mathcal{F}=
\begin{bmatrix}
O_{4\times 4} & \mathcal{F}^1 \\
\\
\mathcal{F}^2 & 0
\end{bmatrix}
,
\end{equation}
where
\begin{subequations}
\begin{align}
\mathcal{F}^1&=
\begin{bmatrix}
\mathcal{L}_{02}^W \\
\\
0
\end{bmatrix}
=
\begin{bmatrix}
\frac{{c_b}}{k_2}{\partial_1}(k_2\cdot) \\
\\
\frac{{c_b}}{k_2}{\partial_2}(k_2\cdot) \\
\\
\frac{{c_b}}{k_2}{\partial_3}(k_2\cdot) \\
\\
0
\end{bmatrix}, \\
\mathcal{F}^2&=
\Big[ -\mathcal{W}\restriction^2 \qquad 0\Big]= \nonumber \\
&=
\begin{bmatrix}
\frac{{c_b}}{k_1}{\partial_1}( k_1 \cdot) &
\frac{{c_b}}{k_1}{\partial_2}(k_1 \cdot) &
\frac{{c_b}}{k_1}{\partial_3}(k_1 \cdot) & 0
\end{bmatrix}.
\end{align}
\end{subequations}\\

The domains of the operators of $\bm{F}^1, \bm{F}^2$, which realize $\mathcal{F}^1, \mathcal{F}^2$, should enjoy
\begin{align}
&\mathsf{D}(\bm{F})=\mathsf{D}(\bm{F}^2)\times \mathsf{D}(\bm{F}^1)= \nonumber \\
&=\mathsf{D}(\bm{C})=\mathfrak{G}_0 \times \mathsf{D}(\bm{L}^W)
=\mathfrak{G}_0\times (\mathfrak{h}^1\times \mathfrak{f})
=(\mathfrak{G}_0\times \mathfrak{h}^1)\times \mathfrak{f}.
\end{align}

Hence we have
\begin{subequations}
\begin{align}
\mathsf{D}(\bm{F}^1)&=\mathfrak{f}=
\Big\{ 
w \in \mathfrak{h}^1 \Big|\quad \mathcal{L}_{02}^Ww=
\frac{c_b}{k_2}\mathrm{grad}(k_2 w) \in \mathfrak{H} \Big\}, \nonumber \\
&\bm{F}^1w=
\begin{bmatrix}
\mathcal{L}_{02}^Ww \\
\\
0
\end{bmatrix}
=\begin{bmatrix}
\frac{c_b}{k_2}\mathrm{grad}(k_2w) \\
\\
0
\end{bmatrix}
\in \mathfrak{H} \times \{0\}
, \\
\mathsf{D}(\bm{F}^2)&=\mathfrak{G}_0\times \mathfrak{h}^1, \nonumber \\
&\bm{F}^2
\begin{bmatrix}
\bm{v} \\
\\
w
\end{bmatrix}
=-\mathcal{W}\restriction^2\bm{v}=\frac{c_b}{k_1}\mathrm{div}(k_1\bm{v})
\in \mathfrak{h}^1
\end{align}
\end{subequations}\\

We claim

\begin{Lemma}
The operator $\mathrm{i}\bm{F}$ is symmetric, that is, $\bm{F}^1 \subset -(\bm{F}^2)^*$ and
$\bm{F}^2 \subset -(\bm{F}^1)^*$. 
\end{Lemma}

Proof. We claim that it holds
$$
\Big(\bm{F}^1 w^2\Big|
\begin{bmatrix}
\bm{v} \\
w^1
\end{bmatrix}
\Big)_{\mathfrak{H}\times \mathfrak{h}^1}
=-\Big(
w^2\Big|\bm{F}^2
\begin{bmatrix}
\bm{v} \\
w^1
\end{bmatrix}
\Big)_{\mathfrak{h}^1}
$$
for $\forall w^2 \in \mathsf{D}(\bm{F}^1), \forall 
\begin{bmatrix}
\bm{v} \\
w^1
\end{bmatrix}
\in \mathsf{D}(\bm{F}^2)$.
But 
\begin{align*}
&\mbox{Left-hand side}=\int_{\mathfrak{R}_b}(\mathrm{grad}(k_2 w^2)| k_1\bm{v})d\bm{x}, \\
&\mbox{Right-hand side}=- \int_{\mathfrak{R}_b}k_2w^2\cdot \mathrm{div}(k_1\bm{v})^*d\bm{x}.
\end{align*}
These are equal since $\bm{v} \in \mathfrak{G}_0$. $\square$\\

Let us look at $\bm{F}^2\bm{F}^1$, an operator in $\mathfrak{h}^1$. By definition we see
\begin{align}
\mathsf{D}(\bm{F}^2\bm{F}^1)&=\Big\{ w \in \mathfrak{h}^1\Big|\quad \frac{c_b}{k_2}\mathrm{grad}(k_2w)
=\frac{k_1}{\rho_b}\mathrm{grad}(k_2w) \in \mathfrak{G}_0
\Big\}, \\
\bm{F}^2\bm{F}^1w&=\frac{c_b}{k_1}\mathrm{div}\Big(\frac{c_bk_1}{k_2}\mathrm{grad}(k_2w)\Big) .
\end{align}

Recall that
$$ \frac{c_b}{k_1}=\frac{1}{\rho_b}\cdot c_b\cdot \frac{1}{E},
\quad \frac{c_bk_1}{k_2}=\rho_b\cdot E^2=\sigma_b\cdot\frac{1}{E^2}, \quad
k_2=c_b\cdot \frac{1}{E}.
$$
Therefore 
$$ \|\bm{v}\|_{\mathfrak{G}}=\Big[ \|\bm{v}\|_{\mathfrak{H}}^2+ \|\mathrm{div}(\rho_b\bm{v})\|_{L^2(\sigma_bd\bm{x})}^2 \Big]^{\frac{1}{2}} $$
is equivalent to
$$
\Big[ 
\|\bm{v}\|_{\mathfrak{H}}^2+ \|\mathrm{div}(\rho_b\bm{v})\|_{L^2(\frac{c_bk_2}{k_1}d\bm{x})}^2 \Big]^{\frac{1}{2}}
$$
for $\bm{v}=\frac{c_b}{k_2}\mathrm{grad}(k_2w)=\frac{k_1}{\rho_b}\mathrm{grad}(k_2w) \in \mathfrak{G}_0$.

We can claim 

\begin{Lemma}
The operator  $-\bm{F}^2\bm{F}^1$ is the self-adjoint operator in 
$\mathfrak{h}^1$
associated with the quadratic form
$$
Q[w]=\int_{\mathfrak{R}_b}\frac{c_bk_1}{k_2}\|\mathrm{grad}(k_2w)\|^2d\bm{x} =(-\bm{F}^1\bm{F}^2w|w)_{\mathfrak{h}^1}.
$$
Moreover the resolvent of $-\bm{F}^2\bm{F}^1$ is compact, therefore the spectrum is of the Sturm-Liouville type, that is, the imbedding $\{ w| \|w\|_{\mathfrak{h}}^2+Q[w] <\infty\} \hookrightarrow \mathfrak{h}^1$ is compact.
\end{Lemma}

Proof. Since
$$
\frac{1}{C}\mathsf{d}^{\frac{1}{\gamma-1}-1}\leq \frac{\rho_b}{(k_2)^2}\leq C \mathsf{d}^{\frac{1}{\gamma-1}-1}, \quad
\frac{1}{C}\leq \Big(\frac{c_b}{k_2}\Big)^2 \leq C,
$$
where $\mathsf{d} =\mathrm{dist}(\cdot, \partial \mathfrak{R})$,
we see that $\|w\|_{\mathfrak{h}}^2+Q[w]$ is equivalent to
$$\|\hat{w}\|_{L^2 (\mathsf{d}^{\frac{1}{\gamma-1}-1}) }^2+
\|\mathrm{grad} \hat{w}\|_{L^2 (\mathsf{d}^{\frac{1}{\gamma-1}}) }^2,$$
where $\hat{w}=k_2w$.
It is known that $W_0^1(\mathsf{d}^{\frac{1}{\gamma-1}-1}, \mathsf{d}^{\frac{1}{\gamma-1}})$ is imbedded compactly into
$L^2(\mathsf{d}^{\frac{1}{\gamma-1}-1})$. ( \cite[Theorem 2.4, or p.740. B]{GurkaO}.))
$\|\hat{w}\|_{L^2(\mathsf{d}^{\frac{1}{\gamma-1}-1})}$ is equivalent to $\|w\|_{\mathfrak{h}}$. $\square$.\\

Let us observe the operator $\bm{F}^1\bm{F}^2$ in $\mathfrak{H}\times \mathfrak{h}^1$:
$$\mathsf{D}(\bm{F}^1\bm{F}^2)=\Big\{ 
\begin{bmatrix}
\bm{v} \\
w
\end{bmatrix}
\Big|\quad \bm{v} \in \mathfrak{G}_0, 
\frac{{c_b}}{k_2}\mathrm{grad}
\Big(\frac{{c_b} k_2}{k_1}\mathrm{div}(k_1 \bm{v})\Big) 
\in \mathfrak{H}
\Big\},
$$

$$\bm{F}^1\bm{F}^2
\begin{bmatrix}
\bm{v} \\w
\end{bmatrix}
=
\begin{bmatrix}
\frac{{c_b}}{k_2}\mathrm{grad}
\Big(\frac{{c_b} k_2}{k_1}\mathrm{div}(k_1 \bm{v})\Big) \\
0
\end{bmatrix}.
$$

We note that, if $\mathrm{div}(k_1\bm{v})=0$, then
$\begin{bmatrix}
\bm{v} \\
w
\end{bmatrix}
\in \mathsf{N}(\bm{F}^1\bm{F}^2)$ for any $w \in 
\mathfrak{h}$.
Therefore the dimension of the null space is infinity.\\

We see 
\begin{equation}
-\bm{F}^1\bm{F}^2=
\begin{bmatrix}
\bm{L}^{\sharp} & \bm{0} \\
\\
\bm{0}^{\perp} & 0 
\end{bmatrix}
,
\end{equation}
where $\bm{L}^{\sharp}$ is the operator in $\mathfrak{H}$ defined as
\begin{align}
&\mathsf{D}(\bm{L}^{\sharp})=\Big\{ \bm{v} \in \mathfrak{G}_0 \Big| \quad \mathcal{L}^{\sharp}\bm{v} \in \mathfrak{H} \Big\}, \\ 
& \bm{L}^{\sharp}\bm{v}=\mathcal{L}^{\sharp}\bm{v}=
-\frac{c_b}{k_2}\mathrm{grad}\Big(\frac{c_bk_2}{k_1}\mathrm{div}(k_1\bm{v})\Big).
\end{align}
Recall that
$$ \frac{c_b}{k_2}=E,\quad \frac{c_bk_2}{k_1} =\sigma_b\cdot \frac{1}{E^2},
\quad k_1=\rho_b\cdot E.
$$
Therefore 
$$ \|\bm{v}\|_{\mathfrak{G}}=\Big[ \|\bm{v}\|_{\mathfrak{H}}^2+ \|\mathrm{div}(\rho_b\bm{v})\|_{L^2(\sigma_bd\bm{x})}^2 \Big]^{\frac{1}{2}} $$
is equivalent to
$$
\Big[ 
\|\bm{v}\|_{\mathfrak{H}}^2+ \|\mathrm{div}(k_1\bm{v})\|_{L^2(\frac{c_bk_2}{k_1}d\bm{x})}^2 \Big]^{\frac{1}{2}}
$$

Thus $\bm{L}^{\sharp}$ is the Friedrichs extension of $\mathcal{L}^{\sharp}\restriction C_0^{\infty}$ associated with the quadratic form
$$
Q^{\sharp}[\bm{v}]=\int_{\mathfrak{R}_b}|\mathrm{div}(k_1\bm{v})|^2\frac{c_bk_2}{k_1}d\bm{x} \quad (\bm{v} \in \mathfrak{G}_0),
$$
for which
$$Q^{\sharp}(\bm{v},\bm{v}')=(\bm{L}^{\sharp}\bm{v}|\bm{v}')_{\mathfrak{H}} \quad (\bm{v} \in \mathsf{D}(\bm{L}^{\sharp}), \bm{v}' \in \mathfrak{G}_0).
$$

Consequently $\bm{L}^{\sharp}$ is a self-adjoint operator in $\mathfrak{H}$ and
$-\bm{F}^1\bm{F}^2$ is a self-adjoint operator in $\mathfrak{H}\times \mathfrak{h}^1$.\\


Summing up, we claim \\

\textbullet\ $\bm{F}^1: (\subset \mathfrak{h}^1) \rightarrow \mathfrak{H}\times \mathfrak{h}^1, \bm{F}^2:(\subset \mathfrak{H} \times \mathfrak{h}^1) \rightarrow \mathfrak{h}^1$ are densely defined;

\textbullet\ $\bm{F}^1 \subset -(\bm{F}^2)^*, \bm{F}^2\subset -(\bm{F}^1)^*$, and $\mathrm{i}\bm{F}$ is symmetric. 

\textbullet\ $\bm{F}^2\bm{F}^1 : (\subset \mathfrak{h}^1) \rightarrow \mathfrak{h}^1$ is a self-adjoint operator in $\mathfrak{h}^1$ and the spectrum is of the Strum-Liouville type, $\mbox{\Pisymbol{psy}{83}}_e(\bm{F}^2\bm{F}^1) =\emptyset$.

\textbullet\ $\bm{F}^1\bm{F}^2$ is a self-adjoint operator in $\mathfrak{H} \times \mathfrak{h}^1$.\\

Here the essential spectrum $\mbox{\Pisymbol{psy}{83}}_e$ is defined as follows:

\begin{Definition}
A densely defined closed operator $T$ in a Banach space $\mathsf{X}$ is said to be Fredholm if both $\mathrm{dim}(\mathsf{N}(T)) $ and 
$\mathrm{dim}(\mathsf{X}/\mathsf{R}(T)) $ are finite. 

The essential spectrum $\mbox{\Pisymbol{psy}{83}}_e(T)$ is the set of complex numbers $\lambda$ such that $T-\lambda$ is not Fredholm.
\end{Definition}

\begin{Remark}
It is known that if $T$ is Fredholm, then $\mathsf{R}(T)$ is closed. (\cite[ Section IV.5.1. p. 230]{Kato}) Therefore the densely defined closed operator $T$ is Fredholm if and only if $\mathsf{R}(T)$ is closed and both
$\mathrm{dim}(\mathsf{N}(T)) $ and 
$\mathrm{dim}(\mathsf{X}/\mathsf{R}(T)) $ are finite. 

It is known that $\mbox{\Pisymbol{psy}{83}}_e(T)=\mbox{\Pisymbol{psy}{83}}_e(T+A)$ for any compcat operator $A$ on $\mathsf{X}$. (\cite[Theorem IV.5.26.]{Kato})
\end{Remark}

Then the M\"{o}ller's theory \cite{Moeller} is applicable to claim:\\

\textbullet\ $\mbox{\Pisymbol{psy}{83}}({\bm{F}^1\bm{F}^2})$ is a discrete subset of $\mathbb{C}$ and $\mbox{\Pisymbol{psy}{83}}_e(
{\bm{F}^1\bm{F}^2}) \subset \{ 0 \}$.

\textbullet\ $\mbox{\Pisymbol{psy}{83}}({\bm{F}})$ is a discrere subset of $\mathbb{C}$, and
$\mbox{\Pisymbol{psy}{83}}_e({\bm{F}}) \subset \{ 0 \}$.

\textbullet\ It holds that
$\mbox{\Pisymbol{psy}{83}}_e({\bm{F}}+\bm{H})=\mbox{\Pisymbol{psy}{83}}_e(\mathcal{J}^*\bm{H}\mathcal{J})$, where $\mathcal{J}$ is the canonical imbedding of $\mathsf{N}({\bm{F}^1\bm{F}^2})$ into $\mathfrak{E}_{W}=\mathfrak{H}\times \mathfrak{h}^2$,  namely 
$$\mathcal{J}
\begin{bmatrix}
\bm{v} \\
\\
w^1
\end{bmatrix}
=
\begin{bmatrix}
\bm{v} \\
\\
w^1 \\
\\
0
\end{bmatrix}
\quad\mbox{for}\quad
\mathrm{grad}\Big(\frac{{c_b}}{k_1}\mathrm{div}(k_1 \bm{v})\Big)=\bm{0}.
$$

Note 

\begin{Lemma}
${\mathcal{J}}=
\begin{bmatrix}
{{\iota}} \\
\\
0
\end{bmatrix}
$ 
, where
${{\iota}} $ is the canonical imbedding of $\mathsf{N}({{\bm{F}}^1{\bm{F}}^2})$
into $\mathfrak{H} \times \mathfrak{h}^1$, and  $ {\mathcal{J}}^*=
\begin{bmatrix}
{{\iota}}^* & 0
\end{bmatrix}
$, while ${{\iota}}^*$ is the orthogonal projection on $\mathfrak{H}\times \mathfrak{h}^1$ onto $\mathsf{N}({{\bm{F}}^1{\bm{F}}^2})$.
\end{Lemma}

\section{Essential spectrum of the operator of {\bf ELASO\ddag($\lozenge$)} }

We consider the essential spectrum of the operator $\bm{C}$. 
Since $\bm{G}$ is compact, we have
$\mbox{\Pisymbol{psy}{83}}_e({\bm{C}})=\mbox{\Pisymbol{psy}{83}}_e({\bm{F}}+\bm{H})$.
Moreover we have verified that $\mbox{\Pisymbol{psy}{83}}_e({\bm{F}}+\bm{H})=\mbox{\Pisymbol{psy}{83}}_e(\mathcal{J}^*\bm{H}\mathcal{J})$, with
the canonical imbedding $\mathcal{J}$ of $\mathsf{N}({\bm{F}^1\bm{F}^2})$ into $\mathfrak{E}_{W}=\mathfrak{H}\times \mathfrak{h}^2$. Therefore we have to
analyze $\mbox{\Pisymbol{psy}{83}}_e(\mathcal{J}^*\bm{H}\mathcal{J})$ .\\

In order to analyze $\mbox{\Pisymbol{psy}{83}}_e(\mathcal{J}^*\bm{H}\mathcal{J})$ 
we make use of the cylindrical co-ordinate $(\varpi, z, \phi)$:
$$x^1=\varpi\cos\phi, \quad x^2=\varpi\sin\phi, \quad x^3=z.
$$ 
The variable $\bm{v}=(v^1, v^2, v^3)^{\top}$ is transformed to the variable
$ \grave{\bm{v}}=(v^{\varpi}, v^z, v^{\phi})^{\top}$ by
\begin{equation}
\bm{v}=\bm{P}\grave{\bm{v}}=
\begin{bmatrix}
v^1 \\
\\
v^2 \\
\\
v^3
\end{bmatrix}
=
\begin{bmatrix}
\cos\phi & 0 & -\sin\phi \\
\\
\sin\phi & 0 & \cos\phi \\
\\
0 & 1 & 0
\end{bmatrix}
\begin{bmatrix}
v^{\varpi} \\
\\
v^{z} \\
\\
v^{\phi}
\end{bmatrix}
\end{equation}
We denote
\begin{equation}
\grave{W}=
\begin{bmatrix}
v^{\varpi} \\
v^z \\
v^{\phi} \\
w^1 \\
w^2
\end{bmatrix}
\end{equation}
and
\begin{align}
\bm{P}_4&=
\begin{bmatrix}
\bm{P} & O \\
\\
0 & 1
\end{bmatrix}:
\mathfrak{h}_{\grave{\bm{v}}, w^1}^4 \rightarrow \mathfrak{h}_{\grave{\bm{v}}, w^1}^4, \\
\bm{P}_5&=
\begin{bmatrix}
\bm{P} & O \\
\\
O & I_{2\times 2}
\end{bmatrix}:
\mathfrak{h}_{\grave{W}}^5 \rightarrow \mathfrak{h}_W^5.
\end{align}

Put
\begin{equation}
\grave{\bm{H}}=\bm{P}_5^{-1}\bm{H}\bm{P}_5.
\end{equation}

Then we have
\begin{equation}
\grave{\mathcal{H}}=
\begin{bmatrix}
0 & 0 & -2\Omega & -\frac{{c_b}^2}{k_1}\frac{\partial k_1}{\partial \varpi} & 0\\
\\
0 & 0 & 0 & -\frac{{c_b}^2}{k_1}\frac{\partial k_1}{\partial z} & 0\\
\\
2\Omega & 0 & 0 & 0 & 0 \\
\\
\frac{1}{k_3}\frac{\partial k_3}{\partial \varpi} & \frac{1}{k_3}\frac{\partial k_3}{\partial z} & 0 & 0 & 0 \\
\\
0 & 0 & 0 & 0 & 0
\end{bmatrix},
\end{equation}
and
\begin{equation}
\mathcal{J}^*\bm{H}\mathcal{J}=\bm{P}_5\grave{\mathcal{J}}^*\grave{\bm{H}}\grave{\mathcal{J}}\bm{P}_5^{-1},
\qquad
\mbox{\Pisymbol{psy}{83}}_e(\mathcal{J}^*\bm{H}\mathcal{J})=
\mbox{\Pisymbol{psy}{83}}_e(\grave{\mathcal{J}}^*\grave{\bm{H}}\grave{\mathcal{J}}),
\end{equation}
where $\grave{\mathcal{J}}=\bm{P}_5^{-1}\mathcal{J}\bm{P}_5$ turns out to be the canonical imbedding of
$\mathsf{N}({\grave{\bm{F}}^1\grave{\bm{F}}^2)}$ into $\mathfrak{h}_{\grave{W}}^5$, namely
$$
\grave{\mathcal{J}}
\begin{bmatrix}
\grave{\bm{v}} \\
\\
w^1
\end{bmatrix}
=
\begin{bmatrix}
\grave{\bm{v}} \\
\\
w^1 \\
\\
0
\end{bmatrix}
\quad\mbox{for}\quad
\mathrm{Grad}\Big(\frac{{c_b}k_2}{k_1}\mathrm{Div}(k_1\grave{\bm{v}})\Big)=\bm{0}.
$$

Here
we denote
\begin{equation}
\mathrm{Grad}q=
\begin{bmatrix}
\frac{\partial q}{\partial\varpi} \\
\\
\frac{\partial q}{\partial z} \\
\\
\frac{1}{\varpi}\frac{\partial q}{\partial\phi}
\end{bmatrix},
\end{equation}
and
\begin{equation}
\mathrm{Div}
\begin{bmatrix}
q^{\varpi} \\
\\
q^z \\
\\
q^{\phi}
\end{bmatrix}
=
\frac{1}{\varpi}\frac{\partial}{\partial \varpi}(\varpi q^{\varpi})+\frac{\partial}{\partial z}q^z
+\frac{1}{\varpi}
\frac{\partial}{\partial \phi}q^\phi.
\end{equation}

Of course, the operators $\grave{\bm{F}}^1: (\subset \mathfrak{h}_{w^2}^1) \rightarrow \mathfrak{h}_{\grave{\bm{v}}, w^1}^4$ and $\grave{\bm{F}}^2: (\subset \mathfrak{h}_{\grave{\bm{v}}, w^1}^4) \rightarrow \mathfrak{h}_{w^2}^1$ are defined as $\bm{F}^1, \bm{F}^2$ by using
$\mathrm{Grad}, \mathrm{Div}$ instead of $\mathrm{grad}, \mathrm{div}$
, namely, $\grave{\bm{F}}^1=\bm{P}_4^{-1}\bm{F}^1$ and $\grave{\bm{F}}^2=\bm{F}^2\bm{P}_4$.
It holds
\begin{equation}
\grave{\bm{F}}=\bm{P}_5^{-1}\bm{F}\bm{P}_5=
\begin{bmatrix}
O & \grave{\bm{F}}^1 \\
\\
\grave{\bm{F}}^2 & 0
\end{bmatrix}
=\begin{bmatrix}
O & \bm{P}_4^{-1}\bm{F}^1 \\
\\
\bm{F}^2\bm{P}_4 & 0
\end{bmatrix},
\end{equation}
while
\begin{subequations}
\begin{align}
& \grave{\bm{F}}^1=
\begin{bmatrix}
\frac{c_b}{k_2}\frac{\partial}{\partial\varpi}(k_2\cdot) \\
\\
\frac{c_b}{k_2}\frac{\partial}{\partial z}(k_2\cdot) \\
\\
\frac{c_b}{k_2\varpi }\frac{\partial}{\partial \phi}(k_2\cdot) \\
\\
0
\end{bmatrix}, \\
&\grave{\bm{F}}^2=
\Big[
\frac{c_b}{k_1\varpi}\frac{\partial}{\partial \varpi}(\varpi k_1\cdot) \quad
\frac{c_b}{k_1}\frac{\partial}{\partial z}( k_1\cdot) \quad
\frac{c_b}{k_1\varpi}\frac{\partial}{\partial \phi}( k_1\cdot) \quad
0
\Big].
\end{align}
\end{subequations}\\

Let us look at 
$$
\grave{\bm{H}}=
\begin{bmatrix}
O_{2\times 2} & \bm{H}^1 & O_{2\times 1} \\
\\
\bm{H}^2 & O_{2\times 2} & O_{2\times 1} \\
\\
O_{1\times 2} & O_{1\times 2} & 0 
\end{bmatrix},
$$
where
$(\bm{H}^j\bm{u})({\bm{x}})=H^j({\bm{x}})\bm{u}({\bm{x}})$ for $\bm{u}\in \mathfrak{h}^2, 
j=1,2$, 
$$
H^1=
\begin{bmatrix}
\displaystyle -2\Omega & \displaystyle -\frac{{c_b}^2}{k_1}\frac{\partial k_1}{\partial \varpi} \\
\\
\displaystyle 0 &\displaystyle -\frac{{c_b}^2}{k_1}\frac{\partial k_1}{\partial z}
\end{bmatrix}
, \quad
H^2=
\begin{bmatrix}
2\Omega & 0 \\
\\
\\
\displaystyle \frac{1}{k_3}\frac{\partial k_3}{\partial \varpi} & \displaystyle \frac{1}{k_3}\frac{\partial k_3}{\partial z}
\end{bmatrix}.
$$

Note that
$\grave{\mathcal{J}}=
\begin{bmatrix}
{\grave{\iota}} \\
\\
0
\end{bmatrix}
$ 
, where
${\grave{\iota}} $ is the canonical imbedding of $\mathsf{N}({\grave{\bm{F}}^1\grave{\bm{F}}^2)}$
into $\mathfrak{h}_{\grave{\bm{v}}, w^1}^4$, and that $ \grave{\mathcal{J}}^*=
\begin{bmatrix}
{\grave{\iota}}^* & 0
\end{bmatrix}
$, while ${\grave{\iota}}^*$ is the orthogonal projection on $\mathfrak{h}_{\grave{\bm{v}}, w^1}^4$ onto $\mathsf{N}(
\overline{\grave{\bm{F}}^1\grave{\bm{F}}^2)}$.
Thus we should consider 
\begin{equation}
\bm{M}:={\grave{\iota}}^*
\begin{bmatrix}
O & \bm{H}^1 \\
\bm{H}^2 & O
\end{bmatrix}
{\grave{\iota}} =
\grave{\mathcal{J}}^*
\grave{\bm{H}}
\grave{\mathcal{J}}
\end{equation}
and study the essential spectrum of the bounded operator $\bm{M}$ on 
$\mathsf{N}({\grave{\bm{F}}^1\grave{\bm{F}}^2)}$,
for which $\mbox{\Pisymbol{psy}{83}}_e(\bm{M})=
\mbox{\Pisymbol{psy}{83}}_e(\grave{\mathcal{J}}^*\grave{\bm{H}}\grave{\mathcal{J}})=
\mbox{\Pisymbol{psy}{83}}_e({\mathcal{J}}^*{\bm{H}}{\mathcal{J}})=
\mbox{\Pisymbol{psy}{83}}_e({\bm{C}})$.\\

We are going to apply the following theorem (see e.g., \cite[Chapter 1, Corollary 4.7 ]{Edmunds}):
\begin{quote}

{\it Let $T$ be a bounded linear operator in a Hilbert space $\mathsf{X}$. If there is a sequence $(x_n)_n$ in $\mathsf{X}$ sich that
$\|x_n\|\geq \frac{1}{C}>0$, $x_n \rightharpoonup 0$ weakly, and
$\|(T-\lambda)x_n\| \rightarrow 0$ as $n \rightarrow \infty$, then 
$T-\lambda$ is not Fredholm, $\lambda \in \mbox{\Pisymbol{psy}{83}}_e(T)$. }
\end{quote}

Such a sequence is called `Weyl's singular sequence'.\\

Let us consider $\lambda \not=0$
and look at
$$\lambda-\bm{M}={\grave{\iota}}^*
\begin{bmatrix}
\lambda & -\bm{H}^1 \\
\\
-\bm{H}^2 & \lambda
\end{bmatrix}
{\grave{\iota}}.
$$ \\

\begin{Definition}
Let us denote by $\alpha_{\pm}(\bm{x})$, $\alpha_-(\bm{x})\leq \alpha_+(\bm{x}),$ the eigenvalues of the symmetric matrix
$$-\frac{1}{2}( H^1H^2+(H^1H^2)^{\top})(\bm{x}).$$
\end{Definition}

More concretely, we put
\begin{equation}
\alpha_{\pm}=\frac{1}{2}\Big( q_1\pm\sqrt{ (q_1)^2+(q_2)^2 } \Big),
\end{equation}
where
\begin{subequations}
\begin{align}
q_1&=4\Omega^2+
\frac{{c_b}^2}{k_1k_3}\Big(\frac{\partial k_1}{\partial \varpi}\frac{\partial k_3}{\partial\varpi}+\frac{\partial k_1}{\partial z}\frac{\partial k_3}{\partial z}\Big), \\
q_2&=
\frac{{c_b}^2}{k_1k_3}\Big(\frac{\partial k_1}{\partial \varpi}\frac{\partial k_3}{\partial z}+\frac{\partial k_1}{\partial z}\frac{\partial k_3}{\partial \varpi}\Big).
\end{align}
\end{subequations}\\

We claim
\begin{Theorem}\label{Th.ESP1}
Let $\lambda\not=0$. Suppose that there is ${\bm{x}}_0=(\varpi_0, 0, z_0) \in \mathfrak{R}_b$ with $\varpi_0 >0$
such that $-\lambda^2 \in [\alpha_-({\bm{x}}_0), \alpha_+({\bm{x}}_0)]$. Then $\lambda \in \mbox{\Pisymbol{psy}{83}}_e(\bm{M})=
\mbox{\Pisymbol{psy}{83}}_e({\bm{C}}) $.
\end{Theorem}

Proof of Theorem \ref{Th.ESP1}. Since $-\lambda^2 \in [\alpha_-({\bm{x}}_0), \alpha_+({\bm{x}}_0)]$, there is a vector $\bm{c}_1 \in \mathbb{R}^2$ such that $\|\bm{c}_1\|=1$ and
$$
-\lambda^2 = -\frac{1}{2}\Big(\bm{c}_1\Big| (H^1H^2+(H^1H^2)^{\top})({\bm{x}}_0)\bm{c}_1\Big).
$$
Let us fix such an $\bm{c}_1$, and put
$$
\bm{c}_2:=J\bm{c}_1, \quad J=
\begin{bmatrix}
0 & -1 \\
\\
1 & 0
\end{bmatrix}.
$$
The $(\bm{c}_1, \bm{c}_2)$ is an orthonormal basis of $\mathbb{R}^2$.
Note that it holds
\begin{equation}
\lambda^2=\Big(\bm{c}_1\Big| H^1H^2({\bm{x}}_0)\bm{c}_1\Big) \label{ESP1}
\end{equation}

Let us define, for $0<\varepsilon \ll1$, the function $u_{\varepsilon}^{\flat}
\in C_0^{\infty}(\mathfrak{R}_b; \mathbb{R})$ of the form
$u_{\varepsilon}^{\flat}(\bm{x})=u_{\varepsilon}(\varpi, z), \varpi=\sqrt{(x^1)^2+(x^2)^2}, z=x^3$ as follows.
Let $\varphi \in C_0^{\infty}(\mathbb{R}; \mathbb{R})$ satisfy
$\mathrm{supp}[\varphi]\subset ]-1, 1[$ and
$\displaystyle \int_{-\infty}^{+\infty}|\varphi(\xi)|^2d\xi =1$. 
Let $0<\nu_1<\nu_2$. We introduce the co-ordinates $\xi^1, \xi^2$ on $\mathbb{R}^2=\{ \overset{=}{\bm{x}} =(\varpi, z)\}$ by putting
$$
\xi^1=\varepsilon^{-\nu_1}(\bm{c}_1|\overset{=}{\bm{x}}-\overset{=}{\bm{x}}_0), \quad
\xi^2=\varepsilon^{-\nu_2}(\bm{c}_2|\overset{=}{\bm{x}}-\overset{=}{\bm{x}}_0),
$$
where $\overset{=}{\bm{x}}_0=(\varpi_0, z_0)$. 
Put
$$u_{\varepsilon}(\overset{=}{\bm{x}})=\varepsilon^{\frac{\nu_2-\nu_1}{2}}\varphi(\xi^1)\varphi(\xi^2).
$$
Then 
\begin{align*}
&\|u_{\varepsilon}\|_{\mathfrak{h}} \leq C\varepsilon^{\nu_2}, \\
& |(\bm{c}_1|\nabla u_{\varepsilon})_{\mathfrak{h}} |\leq C\varepsilon^{\nu_2-\nu_1}, \\
& 0<\frac{1}{C} \leq \|\nabla u_{\varepsilon}\|_{\mathfrak{h}} \leq C\\
&\mathrm{supp}[u_{\varepsilon}] \subset
\{ \overset{=}{\bm{x}} \ |\ \|\overset{=}{\bm{x}}-\overset{=}{\bm{x}}_0\| \leq \sqrt{2}\varepsilon^{\nu_1} \}.
\end{align*}
Here $C$ stands for constants independent of $\varepsilon, 0<\varepsilon \leq \varepsilon_*$
and we denote
$$\nabla u =
\begin{bmatrix}
\displaystyle \frac{\partial u}{\partial\varpi} \\
\\
\displaystyle \frac{\partial u}{\partial z}
\end{bmatrix}.
$$
We are considering small $\varepsilon_*$ such that
$$\{ \overset{=}{\bm{x}} \ |\ \|\overset{=}{\bm{x}}-\overset{=}{\bm{x}}_0\| \leq \sqrt{2}\varepsilon_*^{\nu_1} \} 
\subset \overset{=}{\mathfrak{R}} \cap \{ \varpi >0\}, $$
when $ 0<\frac{1}{C}\leq \rho_b\varpi \leq C$ on $
\{ \overset{=}{\bm{x}} \ |\ \|\overset{=}{\bm{x}}-\overset{=}{\bm{x}}_0\| \leq \sqrt{2}\varepsilon_*^{\nu_1} \} $ so that
$$\frac{1}{C}\|u\|_{L^2} \leq \|u\|_{\mathfrak{h}}
\leq C\|u\|_{L^2} $$
if $\mathrm{supp}[u] \subset \{ \overset{=}{\bm{x}} \ |\ \|\overset{=}{\bm{x}}-\overset{=}{\bm{x}}_0\| \leq \sqrt{2}\varepsilon_*^{\nu_1} \}$, where $\displaystyle \|u\|_{L^2}=\Big[\int_{\overset{=}{\mathfrak{R}}}|u|^2d\overset{=}{\bm{x}} \Big]^{\frac{1}{2}}
= \Big[\int_{\overset{=}{\mathfrak{R}}}|u|^2d\varpi dz \Big]^{\frac{1}{2}}
$ and $\overset{=}{\mathfrak{R}}=\{ (\varpi, z) |\quad (\varpi, 0, z) \in \mathfrak{R}_b \}$.\\
Denoting
$$\tilde{\nabla}u=
\begin{bmatrix}
\displaystyle \frac{1}{\varpi}\frac{\partial}{\partial \varpi}(\varpi u) \\
\\
\displaystyle \frac{\partial u}{\partial z}
\end{bmatrix}
=
\nabla u+\frac{1}{\varpi}
\begin{bmatrix}
u \\
\\
0
\end{bmatrix},
$$
we consider
$$
\bm{g}_{\varepsilon}:=
\begin{bmatrix}
\displaystyle \frac{1}{k_1}\tilde{\nabla}u_{\varepsilon} \\
\\
g_{\varepsilon}^3 \\
\\
g_{\varepsilon}^4
\end{bmatrix},
$$
where
\begin{align}
g_{\varepsilon}^3&=
\frac{1}{\lambda}
\begin{bmatrix}
2\Omega & 0
\end{bmatrix}
\frac{1}{k_1} J\tilde{\nabla}u_{\varepsilon} , \label{ESP2} \\
g_{\varepsilon}^4&=\frac{1}{\lambda}\frac{1}{k_3}(\nabla k_3)^{\top} 
\frac{1}{k_1}J\tilde{\nabla}u_{\varepsilon}. \label{ESP3}
\end{align}
We see
$$ 0<\frac{1}{C} \leq \Big\| \frac{1}{k_1}\tilde{\nabla}u_{\varepsilon} \Big\|_{\mathfrak{h}^2}
\leq 
\|\bm{g}_{\varepsilon}\|_{\mathfrak{h}^4} $$
and $\bm{g}_{\varepsilon} \rightharpoonup 0$ as $\varepsilon \rightarrow 0$ weakly in $\mathfrak{h}^4$. 
We note that $\bm{g}_{\varepsilon} \in \mathsf{N}(\grave{\bm{F}}^2) \subset \mathsf{N}({\grave{\bm{F}}^1\grave{\bm{F}}^2)}$. In fact, since 
$ \tilde{\nabla}^{\top}J\tilde{\nabla}u_{\varepsilon}=0 $,
we see
$$\bm{F}^2\bm{g}_{\varepsilon}=
\frac{{c_b}}{k_1} \tilde{\nabla}^{\top}J\tilde{\nabla}u_{\varepsilon}
+\frac{{c_b}}{\varpi}\frac{\partial g_{\varepsilon}^3}{\partial \phi}=0 $$
We claim that 
\begin{equation}
(\lambda -\bm{M})\bm{g}_{\varepsilon}=
{\grave{\iota}}^*
\begin{bmatrix}
\lambda & -\bm{H}^1 \\
\\
-\bm{H}^2 & \lambda
\end{bmatrix}
{\grave{\iota}} \bm{g}_{\varepsilon} \rightarrow 0 \quad
\mbox{as}\quad \varepsilon \rightarrow 0 \label{ESP*1}
\end{equation}
in $\mathfrak{h}^4$-norm in $\mathsf{N}({\grave{\bm{F}}^1\grave{\bm{F}}^2)}$.
Then $(\bm{g}_{\varepsilon})_{\varepsilon}$ is a Weyl's singular sequence, whose existence implies that
$\lambda-\bm{M}={\grave{\iota}}^*
\begin{bmatrix}
\lambda & -\bm{H}^1 \\
\\
-\bm{H}^2 & \lambda
\end{bmatrix}
{\grave{\iota}} $ is not Fredholm, and $\lambda \in 
\mbox{\Pisymbol{psy}{83}}_e(\bm{M})$. \\
Let $a$ be a real parameter specified later.
The condition
\begin{equation}
\begin{bmatrix}
\lambda & -\bm{H}^1 \\
\\
-\bm{H}^2 & \lambda
\end{bmatrix}
\bm{g}_{\varepsilon}+a \grave{\bm{F}}^1\frac{u_{\varepsilon}}{{c_b}} =o(1) \label{ESP*2}
\end{equation}
implies \eqref{ESP*1}, since ${\grave{\iota}}^*\grave{\bm{F}}^1=0$.
Multiplying by 
$
\begin{bmatrix}
\lambda & \bm{H}^1 \\
\\
O & I
\end{bmatrix}
$ from the left, we see that 
\eqref{ESP*2} is equivalent to
\begin{equation}
\begin{bmatrix}
\lambda^2 -\bm{H}^1\bm{H}^2 & O \\
\\
-\bm{H}^2 & \lambda
\end{bmatrix}
\bm{g}_{\varepsilon}+
a
\begin{bmatrix}
\lambda & \bm{H}^1 \\
\\
O & I
\end{bmatrix}
\grave{\bm{F}}^1\frac{u_{\varepsilon}}{{c_b}} =o(1). \label{ESP*3}
\end{equation}
The 1st and 2nd components of \eqref{ESP*3} read
\begin{align}
&\frac{1}{k_1}(\lambda^2-H^1H^2)J\tilde{\nabla} u_{\varepsilon} +
\frac{a\lambda{c_b}}{k_2}\nabla\Big(\frac{k_2}{{c_b}}u_{\varepsilon}\Big) \nonumber \\
&=\frac{1}{k_1}(\lambda^2-H^1H^2)J\nabla u_{\varepsilon} + a\lambda \nabla u_{\varepsilon} +o(1) \nonumber \\
&=\Big(\frac{1}{k_1}(\lambda^2-H^1H^2)J+a\lambda\Big)\bm{c}_2(\bm{c}_2|\nabla u_{\varepsilon})+o(1)=o(1), \label{ESP**1}
\end{align}
since
$$
\nabla u_{\varepsilon}=\bm{c}_1(\bm{c}_1|\nabla u_{\varepsilon})+
\bm{c}_2(\bm{c}_2|\nabla u_{\varepsilon})=
O(\varepsilon^{\nu_2-\nu_1})+
\bm{c}_2(\bm{c}_2|\nabla u_{\varepsilon}).$$
But \eqref{ESP**1} holds if 
we specify $a$ so that
\begin{equation}
\frac{1}{k_1}(\lambda^2 -H^1H^2)J\bm{c}_2=-a\lambda\bm{c}_2
\quad\mbox{at}\quad \overset{=}{\bm{x}}=\overset{=}{\bm{x}}_0, \label{ESP5}
\end{equation}
which is possible thanks to \eqref{ESP1}. In fact
$$\Big(\frac{1}{k_1}(\lambda^2-H^1H^2)J+a\lambda\Big)\bm{c}_2
=
\Big(\frac{1}{k_1}(\lambda^2-H^1H^2)J+a\lambda\Big)\bm{c}_2\Big|_{\overset{=}{\bm{x}}=\overset{=}{\bm{x}}_0}
+o(1),
$$
since $\mathrm{supp}[u_{\varepsilon}] \rightarrow \{\overset{=}{\bm{x}}_0\}$. 
The 3rd component of \eqref{ESP*3} reads
$$
-
\begin{bmatrix}
2\Omega & 0
\end{bmatrix}
\frac{1}{k_1}J\tilde{\nabla}u_{\varepsilon}+\lambda g_{\varepsilon}^3=o(1).
$$
This condition holds by $g_{\varepsilon}^3$ determined by \eqref{ESP2}.
The 4th component of \eqref{ESP*3} reads
$$
-\frac{1}{k_3}(\nabla k_3)^{\top} J\tilde{\nabla}u_{\varepsilon}+\lambda g_{\varepsilon}^4=o(1).
$$
This condition holds by $g_{\varepsilon}^4$ determined by \eqref{ESP3}.
Summing up, we can claim that \eqref{ESP*3} holds. Therefore we can claim 
\eqref{ESP*1}. 
Proof of Theorem \ref{Th.ESP1} has been completed.\\
\begin{Remark}
The trick of the above discussion is due to
\cite{FLMM} and \cite{FM2000}.
We have followed but little bit simplified 
their settings and proofs.
In fact, the singular sequence $(\bm{g}_{\varepsilon})_{\varepsilon}$ constructed here is such that $\displaystyle \frac{\partial}{\partial \phi}\bm{g}_{\varepsilon}=0$, namely, $\bm{g}_{\varepsilon}(\bm{x})=\bm{g}_{\varepsilon}^{\sharp}(\varpi, z)$. But a singular sequence of the form
$$\bm{g}_{\varepsilon}(\bm{x})=\check{\bm{g}}_{\varepsilon}^{\sharp}(\varpi, z)e^{\mathrm{i}m\phi}$$
with $\bm{x}=(\varpi\cos\phi, \varpi\sin\phi, z)$ can be constructed for the azimuthal wave number $m \in \mathbb{Z}\setminus \{0\}$.
This can be done by taking
$$
\check{\bm{g}}_{\varepsilon}^{\sharp}=
\begin{bmatrix}
\bm{f}_{\varepsilon} \\
\\
\displaystyle \frac{\mathrm{i}\varpi}{m}\frac{1}{k_1}\tilde{\nabla}^{\top} k_1\bm{f}_{\varepsilon} \\
\\
\displaystyle \frac{1}{\lambda}\frac{1}{k_3}(\nabla k_3)^{\top}\bm{f}_{\varepsilon}
\end{bmatrix}
,
$$
where
$$
\bm{f}_{\varepsilon}=\frac{1}{k_1}J\tilde{\nabla}u_{\varepsilon} -\frac{m}{\mathrm{i}\lambda \varpi_0 k_1(\bm{x}_0)}
(
\begin{bmatrix}
2\Omega & 0
\end{bmatrix}
\bm{c}_1)
u_{\varepsilon} \bm{c}_2, 
$$
\end{Remark}

\begin{Corollary}
Suppose $\Omega=0$ and $\mathfrak{a}_b \not=0$, so that
$\mathfrak{a}_b(\bm{x}_0)\not=0$, or, $\displaystyle \frac{dS_b}{dr}\not=0$ at some $\bm{x}_0 \in 
\mathfrak{R}_b\setminus \{O\}$. Then 
$\alpha_-(\bm{x}_0) <0 <\alpha_+(\bm{x}_0)$ and the cross
$$
\mathrm{K}:=\Big[-\sqrt{|\alpha_-({\bm{x}}_0)|}, \sqrt{|\alpha_-({\bm{x}}_0)|}\Big] \cup
\Big[-\sqrt{\alpha_+({\bm{x}}_0)}, \sqrt{\alpha_+({\bm{x}}_0)}\Big]\mathrm{i}$$
is a subset of $ \mbox{\Pisymbol{psy}{83}}_e(\bm{C})$.
\end{Corollary}

Proof.  Recall
$$\alpha_{\pm}=\frac{1}{2}\Big( q_1\pm\sqrt{ (q_1)^2+(q_2)^2 } \Big),
$$
where
\begin{align*}
q_1&=4\Omega^2+
\frac{{c_b}^2}{k_1k_3}\Big(\frac{\partial k_1}{\partial \varpi}\frac{\partial k_3}{\partial\varpi}+\frac{\partial k_1}{\partial z}\frac{\partial k_3}{\partial z}\Big), \\
q_2&=
\frac{{c_b}^2}{k_1k_3}\Big(\frac{\partial k_1}{\partial \varpi}\frac{\partial k_3}{\partial z}+\frac{\partial k_1}{\partial z}\frac{\partial k_3}{\partial \varpi}\Big).
\end{align*} 
Therefore $\alpha_-=\alpha_+$ if and only if $q_1=q_2=0$, and, then $\alpha_{\pm}=0$
and $\{ \lambda | \lambda \not=0, -\lambda^2 \in [\alpha_-,\alpha_+]\}=\emptyset$.
Otherwise $\alpha_-<0<\alpha_+$.
We are supposing that $\Omega=0$ and the background is spherically symmetric.
Since we are supposing that the EOS is $P=\rho^{\gamma}\exp(S/\mathsf{C}_V)$,
we have
\begin{align*}
q_1&=-\frac{1}{\mathsf{C}_V\gamma \rho_b}
\frac{d P_b}{dr}\frac{d S_b}{dr}, \\
q_2&=-\frac{1}{\mathsf{C}_V\gamma \rho_b}
\frac{2\varpi z}{r^2}
\frac{d P_b}{dr}\frac{d S_b}{dr}.
\end{align*}
Since $d P_b/dr <0$, we can claim that $q_1=q_2=0$ everywhere if and only if
$d S_b/dr =0$ everywhere, that is, the background is isentropic. If the background is not isentropic, there is ${\bm{x}}_0 \in {\mathfrak{R}_b}\setminus \{\varpi=0\}$ such that 
$\displaystyle \frac{dS_b}{dr}\big|_{{\bm{x}}_0}\not=0$.Then 
$\alpha_-({\bm{x}}_0) < 0 < \alpha_+({\bm{x}}_0)$ and the set
$\{ \lambda | -\lambda^2 \in [\alpha_-({\bm{x}}_0),\alpha_+({\bm{x}}_0)]\} $ turns out to be the cross
$$
\mathrm{K}=\Big[-\sqrt{|\alpha_-({\bm{x}}_0)|}, \sqrt{|\alpha_-({\bm{x}}_0)|}\Big] \cup
\Big[-\sqrt{\alpha_+({\bm{x}}_0)}, \sqrt{\alpha_+({\bm{x}}_0)}\Big]\mathrm{i},$$
that is $\mathrm{K} \subset \mbox{\Pisymbol{psy}{83}}_e(\bm{C})$.
$\square$ \\

 Let $\lambda \in \mbox{\Pisymbol{psy}{83}}_e(\bm{M}) \setminus \{0\}$. Since $\lambda I-\bm{M}$ is not Fredholm, it holds that $|\lambda|\leq |\|\bm{M}\||_{\mathcal{B}(\mathsf{N}(\bm{F}^1\bm{F}^2))}$, for, otherwise, say, if $\displaystyle \frac{1}{|\lambda|} |\|\bm{M}\||_{\mathcal{B}(\mathsf{N}(\bm{F}^1\bm{F}^2))} < 1$, then
$\lambda - \bm{M}=\displaystyle \lambda (I-\frac{1}{\lambda}\bm{M})$ would have the bounded inverse, and would be Fredholm. Consequently we can claim that
the essential spectrum of $\bm{M}$, $\mbox{\Pisymbol{psy}{83}}_e(\bm{M})$ is included in the disk
$$
\Big\{ \lambda \in \mathbb{C} \Big| \quad |\lambda|^2 \leq
\Big(4\Omega^2+\Big\|\frac{\nabla P_b}{\rho_b}\Big\|_{L^{\infty}(\mathfrak{R}_b)}^2\Big) \vee
\Big(4\Omega^2+\Big\| \mathfrak{a}_b\Big\|_{L^{\infty}(\mathfrak{R}_b)}^2\Big)
\Big\},
$$
that is, $\mbox{\Pisymbol{psy}{83}}_e(\bm{C})$ is bounded in the $\mathbb{C}$-plane.


\vspace{15mm}

{\bf\Large Acknowledgment}\\

This work was partially done during the stay of the author at Institute of Mathematics, Academia Sinica, in December 2024. The author expresses his sincere thanks to Professor Shih-Hsien Yu for his hospitality and discussion on the problem. 
The author expresses his sincere thanks to the anonimous reviewers for patient and careful reading of the manuscript, indicating errors, and suggesting ameliorations of presentation.
This work is supported by the Research Institute for Mathematical Sciences, International Joint Usage/Research Center located in Kyoto University.\\

{\bf\Large The data availability statement}\\

No new data were created or analyzed in this study.

\vspace{15mm}


\end{document}